\documentclass[a4paper, 12pt]{amsart}
\usepackage{amsmath, amssymb, eucal, amscd, amstext, enumerate, mathrsfs, yfonts, amsfonts}

\newtheorem{thm}{Theorem}[section]
\newtheorem{lem}[thm]{Lemma}

\newtheorem{prop}[thm]{Proposition}
\newtheorem{cor}[thm]{Corollary}
\newtheorem{conjecture}[thm]{Conjecture}

\theoremstyle{definition}
\newtheorem{defn}[thm]{Definition}

\theoremstyle{remark}
\newtheorem{eg}[thm]{Example}
\newtheorem{rem}[thm]{Remark}
\newtheorem{ntn}[thm]{Notation}

\newenvironment{prf}{{\noindent \textbf{Proof:}\ }}{\hfill $\Box$\\ \smallskip}

\numberwithin{equation}{section}

\newcommand{\id}{{\rm id}}
\newcommand{\ti}{\tilde}

\newcommand{\smnoind}{\smallskip\noindent}
\newcommand{\CL}{\mathcal{B}}
\newcommand{\CK}{\mathcal{K}}
\newcommand{\CC}{\mathcal{C}}
\newcommand{\CF}{\mathcal{F}}
\newcommand{\proj}{{\rm Proj}}
\newcommand{\op}{\operatorname{OP}}
\newcommand{\Z}{\operatorname{Z}}
\newcommand{\ta}{\mathfrak{A}}
\newcommand{\tb}{\mathfrak{B}}
\newcommand{\tc}{\mathfrak{C}}
\newcommand{\SF}{\mathfrak{sf}}
\newcommand{\CJ}{\mathcal{J}}

\newcommand{\s}{{\operatorname{sp}}}
\newcommand{\Mv}{{\operatorname{Mv}}}
\newcommand{\her}{\operatorname{her}}
\newcommand{\PZ}{{\operatorname{PZ}}}

\begin{document}

\title{A Murray-von Neumann type classification of $C^*$-algebras}

\author{Chi-Keung Ng  \and Ngai-Ching Wong }

\address[Chi-Keung Ng]{Chern Institute of Mathematics and LPMC, Nankai University, Tianjin 300071, China.}
\email{ckng@nankai.edu.cn}

\address[Ngai-Ching Wong]{Department of Applied Mathematics, National Sun Yat-sen University,
Kaohsiung, 80424, Taiwan.}
\email{wong@math.nsysu.edu.tw}

\thanks{The authors are supported by National Natural Science Foundation of China (11071126), and
Taiwan NSC grant (99-2115-M-110-007-MY3).}

\date{\today}



\begin{abstract}
We define type $\ta$, type $\tb$, type
$\tc$ as well as $C^*$-semi-finite $C^*$-algebras.

It is shown that a von Neumann algebra is a type $\ta$, type $\tb$, type
$\tc$ or $C^*$-semi-finite $C^*$-algebra if and only if it is, respectively, a type I, type II, type III or
semi-finite von Neumann algebra.
Any type I $C^*$-algebra
is of type $\ta$ (actually, type $\ta$ coincides with the discreteness as
defined by Peligrad and Zsid\'{o}), and any type II $C^*$-algebra (as
defined by Cuntz and Pedersen) is of type $\tb$.  Moreover,   any type $\tc$
$C^*$-algebra is of type III (in the sense of Cuntz and Pedersen).  Conversely, any separable purely infinite $C^*$-algebra (in the sense of Kirchberg and R{\o}rdam) with either real rank zero or stable rank one is of type $\tc$.

We also prove that type $\ta$, type $\tb$, type $\tc$ and
$C^*$-semi-finiteness are stable under taking hereditary $C^*$-subalgebras,
multiplier algebras and strong Morita equivalence. Furthermore, any
$C^*$-algebra $A$ contains a largest type $\ta$ closed ideal
$J_\ta$, a largest type $\tb$ closed ideal $J_\tb$, a largest type
$\tc$ closed ideal $J_\tc$ as well as a largest $C^*$-semi-finite
closed ideal $J_\SF$.
Among them, we have $J_\ta + J_\tb$ being an
essential ideal of $J_\SF$, and $J_\ta + J_\tb + J_\tc$ being an
essential ideal of $A$. On the other hand, $A/J_\tc$ is always
$C^*$-semi-finite, and if $A$ is $C^*$-semi-finite, then $A/J_\tb$
is of type $\ta$.
\end{abstract}
\maketitle

\begin{center}
\emph{This paper is dedicated to Charles Batty on the occasion of his 60th birthday.} 
\end{center}

\section{Introduction}

\medskip

In their seminal works (\cite{Murray-vonNeumann36}, see also \cite{Murray90}), Murray and von
Neumann defined three types of  von Neumann algebras (namely, type
I, type II and type III) according to the properties of their
projections.
They showed that any von Neumann algebra is a sum of a type
I, a type II, and a type III
von Neumann subalgebras.
This classification was shown to be very
important and becomes the basic theory for the study of von Neumann algebras (see, e.g., \cite{KRII}).
Since a $C^*$-algebra needs not have any projection, a similar
classification for $C^*$-algebras seems impossible. There is, however, an
interesting classification scheme for $C^*$-algebras proposed by
Cuntz and Pedersen in \cite{CP79}, which  captures some features of
the classification of Murray and von Neumann.

\medskip

The classification theme of $C^*$-algebras took a drastic turn after
an exciting work of Elliott on the classification of $AF$-algebras
through the ordered $K$-theory, in the sense that two $AF$-algebras are isomorphic if and only if
they have the same ordered $K$-theory (\cite{Elliott76}).
Elliott then proposed an invariant consisting of the tracial state space and some $K$-theory datum of the underlying $C^*$-algebra (called the \emph{Elliott invariant}) which could be a suitable candidate for a complete invariant for simple separable nuclear $C^*$-algebras.
Although it is known recently that it is not the case (see \cite{Toms08}), this Elliott invariant still works for a very large class of such $C^*$-algebras (namely, those satisfying certain regularity conditions as described in \cite{ET}).
Many people are still making progress in this direction in trying to find the biggest class of $C^*$-algebras that can
be classified through the Elliott invariant (see, e.g., \cite{Elliott93,Rordam02}).
Notice that this classification is very different from the classification in the sense of Murray and von Neumann.

\medskip

In this article, we reconsider the classification of $C^*$-algebras
through the idea of Murray and von Neumann. Instead of considering
projections in a $C^*$-algebra $A$, we consider open projections and we twist the definition of the finiteness of projections slightly to obtain our classification
scheme.

\medskip

The notion of open projections was introduced by
Akemann (in \cite{Akemann69}).
A projection $p$ in the universal enveloping von Neumann algebra  (i.e.\ the biduals) $A^{**}$ of a $C^*$-algebra $A$ (see, e.g., \cite[\S III.2]{Tak}) is an \emph{open} projection of $A$ if
there is an increasing net $\{a_i\}_{i\in \mathfrak{I}}$ of positive elements in
$A_+$ with $\lim_i a_i = p$ in the $\sigma(A^{**},A^*)$-topology.
In the case when $A$ is commutative, open projections of $A$
are exactly characteristic functions of
open subsets of the spectrum of $A$.
In general, there is a bijective correspondence between open projections of $A$ and hereditary $C^*$-subalgebras of $A$
(where a hereditary $C^*$-subalgebra $B$ corresponds to an open projection $p$ such that $B = pA^{**}p\cap A$; see, e.g., \cite[3.11.10]{Ped79}).
Characterisations and further developments of
open projections can be
found in, e.g.,  \cite{Akemann70, AE02, AP73, Brown88-CJM, Effro63, Ped72, Prosser63}.
Since every element in a $C^*$-algebra is in the closed linear span of
its open projections, it is reasonable to believe that the study of open
projections will provide fruitful information about the underlying
$C^*$-algebra.
Moreover, because of the correspondence between open projections (respectively, central open projections) and hereditary $C^*$-subalgebras (respectively, closed ideals), the notion of strong Morita equivalence as defined by Rieffel (see \cite{Rieffel74} and also \cite{Brown-Rieffel77, Rief}) is found to be very useful in this
scheme.

\medskip

One might wonder why we do not consider the classification of the
universal enveloping von Neumann algebras of $C^*$-algebras to obtain a classification of
$C^*$-algebras. A reason is that for a $C^*$-algebra $A$, its bidual
$A^{**}$ always contains many minimum projections (see, e.g.,
\cite[II.17]{Akemann69}), and hence a reasonable theory of type
classification cannot be obtained without serious modifications.
Furthermore, $A^{**}$ are usually very far away from $A$, and
information of $A$ might not always be respected very well in
$A^{**}$; for example, $c$ and $c_0$ have isomorphic biduals, but
the structure of their open projections can be used to distinguish
them (see, e.g., Example \ref{eg:c-c0} and also Proposition
\ref{prop:pre-op-proj}(b)).

\medskip

As in the case of von Neumann algebras, in order to give a
classification of $C^*$-algebras, one needs, first of all, to
consider a good equivalence relation among   open projections.
After some thoughts and considerations, we end up with the ``spatial
equivalence''  as defined in Section 2, which is weaker than the one
defined by Peligrad and
Zsid\'{o} in \cite{PZ00} and stronger than the ordinary Murray-von
Neumann equivalence.
One reason for making this choice
is that it is precisely the ``hereditarily stable version of
Murray-von Neumann equivalence''  that one might want (see
Proposition \ref{prop:cp-equiv}(a)(5)), and it also coincides with
the ``spatial isomorphism'' of the hereditary $C^*$-subalgebras (see
Proposition \ref{prop:cp-equiv}(a)(2)).

\medskip

Using the spatial equivalence relation, we introduce in Section 3,
the notion of $C^*$-finite $C^*$-algebras.
It is shown that the sum
of all $C^*$-finite hereditary $C^*$-subalgebra is a (not necessarily closed) ideal of
the given $C^*$-algebra. In the case when the $C^*$-algebra is
$\CL(H)$ or $\CK(H)$, this ideal is the ideal of all finite
rank operators on $H$. Moreover, through $C^*$-finiteness, we define
type $\ta$, type $\tb$, type $\tc$ as well as $C^*$-semi-finite
$C^*$-algebras, and we study some properties of them. In particular,
we will show that these properties are stable under taking
hereditary $C^*$-subalgebras, multiplier algebras, unitalization (if
the algebra is not unital) as well as strong Morita equivalence. We
will also show that the notion of  type $\ta$ coincides precisely with the
discreteness as defined in \cite{PZ00}.

\medskip

In Section 4, we will compare these notions with some results in the literature and give some examples.
In particular, we show that
any type I $C^*$-algebra (see, e.g., \cite{Ped79}) is of
type $\ta$; any type II $C^*$-algebra (as defined by Cuntz and Pedersen) is of type $\tb$; any semi-finite $C^*$-algebras (in the sense of Cuntz and Pedersen) is $C^*$-semi-finite; any
purely infinite $C^*$-algebra (in the sense of Kirchberg and R\o rdam) with real rank zero
and any separable purely infinite $C^*$-algebra with stable rank one are of type $\tc$; and any type $\tc$ $C^*$-algebra is of type III (as introduced by Cuntz and Pedersen).
Using our arguments for these results, we also show that any purely infinite $C^*$-algebra is of type III.
Moreover, a von Neumann algebra $M$ is a type $\ta$, a type $\tb$, a type $\tc$
or a $C^*$-semi-finite $C^*$-algebra if and only if $M$ is, respectively, a type I, a type II, a type
III, or a semi-finite von Neumann algebra.

\medskip

In Section 5, we show that any $C^*$-algebra $A$ contains a largest
type $\ta$ closed ideal $J^A_\ta$, a largest type $\tb$ closed ideal
$J^A_\tb$, a largest type $\tc$ closed ideal $J^A_\tc$ as well as a
largest $C^*$-semi-finite closed ideal $J^A_\SF$. It is further shown that
$J^A_\ta + J^A_\tb$ is an essential ideal of $J^A_\SF$, and $J^A_\ta
+ J^A_\tb + J^A_\tc$ is an essential ideal of $A$. On the other
hand, $A/J^A_\tc$ is always a $C^*$-semi-finite $C^*$-algebra, while
$B/J^{B}_\tb$ is always of type $\ta$ if one sets $B:=A/J^A_\tc$.
We also compare $J^{M(A)}_\ta$, $J^{M(A)}_\tb$, $J^{M(A)}_\tc$ and
$J^{M(A)}_\SF$ with $J^{A}_\ta$, $J^{A}_\tb$, $J^{A}_\tc$ and
$J^{A}_\SF$, respectively.

\medskip

\begin{ntn}
Throughout this paper, $A$ is a non-zero $C^*$-algebra, $M(A)$ is the multiplier algebra of $A$, $\Z(A)$ is
the center of $A$, and $A^{**}$ is the bidual of $A$. Furthermore,
$\proj(A)$ is the set of all projections in $A$, while
$\op(A)\subseteq \proj(A^{**})$ is the set of all open projections
of $A$.
All ideals in this paper are two-sided ideals (not assumed to be closed unless specified).

If $x,y\in A^{**}$ and $E$ is a subspace of $A^{**}$, we set $xEy := \{xzy:z\in E\}$, and denote by $\overline{E}$ the norm closure of $E$.
For any $x\in A^{**}$, we set $\her_A(x)$ to be the hereditary $C^*$-subalgebra
$\overline{x^*A^{**}x}\cap A$ of $A$ (note that if $u\in A^{**}$ is a partial isometry, then $\her_A(u) = u^*A^{**}u \cap A = \{x\in A: x=u^*uxu^*u\} = \her_A(u^*u)$).
When $A$ is understood, we will use the notation $\her(x)$ instead.
Moreover, $p_x$ is the right support projection of a
norm one element $x\in A$, i.e.\ $p_x$ is the $\sigma(A^{**},A^*)$-limit of $\{(x^*x)^{1/n}\}_{n\in \mathbb{N}}$ and is the
smallest open projection in $A^{**}$ with $xp_x = x$.
%, then $\her(x)=\her(p_x)$ (see e.g.\ Example \ref{eg:herx-hera} for more information).
\end{ntn}

\medskip

\noindent \emph{Acknowledgement:}
The authors would like to thank L. Brown, E. Effros and G. Elliott
for giving some comments.

\medskip

\section{Spatial equivalence of open projections}

\medskip

In this section, we will consider a suitable equivalence relation on the set of open projections of a $C^*$-algebra.
Let us start with the following example, which shows that the structure of open projections is rich enough to distinguish $c$ and $c_0$, while they have isomorphic biduals (see Proposition \ref{prop:pre-op-proj}(b) below for a more general result).

\medskip

\begin{eg}\label{eg:c-c0}
The sets of open projections of $c_0$ and $c$ can be regarded as the
collections $\mathcal{X}$ and $\mathcal{Y}$, of open subsets of
$\mathbb{N}$ and of open subsets of the one point compactification of
$\mathbb{N}$, respectively.
As ordered sets, $\mathcal{X}$ and
$\mathcal{Y}$ are not isomorphic.
In fact, suppose on the contrary that there
is an order isomorphism $\Psi: \mathcal{Y} \to \mathcal{X}$.
Then $\Psi(\mathbb{N})$ is a proper open
subset of $\mathbb{N}$.
Let $k\notin \Psi(\mathbb{N})$ and $U\in \mathcal{Y}$
with $\Psi(U) = \{k\}$.
As $U$ is a minimal element, it is a singleton set.
Thus, $U\subseteq \mathbb{N}$, which gives the contradiction that $\{k\} \subseteq
\Psi(\mathbb{N})$.
\end{eg}

\medskip

Secondly, we give the following well-known remarks which says that open projections and the
hereditary $C^*$-subalgebras they define, are ``hereditarily
invariant''.
These will clarify some discussions later on.

\medskip

\begin{rem}\label{rem:her-subalg}
Let $B\subseteq A$ be a hereditary $C^*$-subalgebra and $e\in
\op(A)$ be the open projection with $\her_A(e) = B$.

\smnoind (a) For any $p\in \proj(B^{**})$, one has
$\her_B(p) = \her_A(p)$.
%In fact, if $a\in A$ with $a = pap$, then $a = eae$ as $p\leq e$, which implies that $a\in B$. The converse is clear.

\smnoind (b) $\op(B) = \op(A)\cap B^{**}$.
In fact, if $p\in
\op(A)\cap B^{**}$ and $\{a_i\}_{i\in \mathfrak{I}}$ is an
approximate unit in $\her_A(p) = \her_B(p)$, then $\{a_i\}_{i\in \mathfrak{I}}$ will
$\sigma(B^{**}, B^*)$-converge to $p$ and $p\in \op(B)$.
%The converse is obvious.

\smnoind (c) If $z\in A$ satisfying $zz^*,z^*z\in B$, then $z\in B$.
In fact, as $z^*z\in \her_A(e) = eA^{**}e \cap A$, by considering the polar
decomposition of $z$, we see that $ze = z$. Similarly, we have $ez =
z$.

\smnoind
(d) If $f\in \op(A)$, the open projections corresponding to $\her(e)\cap \her(f)$ and the hereditary $C^*$-subalgebra generated by $\her(e) + \her (f)$ are $e\wedge f$ and $e\vee f$ respectively.
\end{rem}

\medskip

Let $j_A:M(A)\to A^{**}$ be the canonical $^*$-monomorphism, i.e.\
$j_A(x)(f) = \ti f(x)$ ($x\in M(A),
f\in A^*$), where $\ti f\in M(A)^*$ is the unique strictly
continuous extension of $f$.
%On the other hand, we denote by $i_A: A\to M(A)$ the canonical embedding.
The proposition below can be regarded as a motivation behind the
study of $C^*$-algebras  through their open projections. It could be
a known result (especially, part (a)).
However, since we need it for the equivalence of (1) and (5) in Proposition \ref{prop:cp-equiv}(a), we give a proof here for completeness.

\medskip

\begin{prop}\label{prop:pre-op-proj}
Suppose that $A$ and $B$ are $C^*$-algebras, and $\Phi: A^{**}\to B^{**}$ is  a $^*$-isomorphism.

\smnoind
(a) If $\Phi\big(j_A(M(A))\big) = j_B(M(B))$, then $\Phi(A) = B$.

\smnoind
(b) If $\Phi(\op(A)) = \op(B)$, then $\Phi(A) = B$.
\end{prop}
\begin{prf}
(a) Let $p_A\in \op(M(A))$ such that $\her_{M(A)}(p_A) = A$.
It is not hard to verify that $p_A$ is the support of $\ti j_A$,
where $\ti j_A: M(A)^{**} \to A^{**}$ is the $^*$-epimorphism induced by $j_A$.
%(this fact also follows from \cite[Theorem 2.2]{Kusuda90}).
%In fact, if $q$ is the support of $\ti j_A$, then $qM(A)^*$ is the set of all strictly continuous functionals. Moreover, $(1-p_A)M(A)^*$ is the kernel of $i_A^*$ (thus, $i_A^*: p_AM(A)^*\to A^*$ is bijective). It is easy to see that $p_A\geq q$ and $i_A^*: qM()A^*\to A^*$ is surjective. Consequently, $p_A = q$.
Consider $\Psi := j_B^{-1}\circ \Phi_{\mid j_A(M(A))}\circ j_A: M(A)\to M(B)$ (which is well-defined by the hypothesis).
Since $j_B\circ \Psi  = \Phi_{\mid j_A(M(A))}\circ j_A$,
we see that $\ti j_B\circ \Psi^{**} = \Phi\circ \ti j_A$ (as $\Phi$ is automatically weak-*-continuous).
Thus, $\ti j_B(\Psi^{**}(p_A)) = 1_{B^{**}}$ which implies $\Psi^{**}(p_A)\geq p_B$.
Similarly,
$$(\Psi^{**})^{-1}(p_B)
\ =\ (j_A^{-1}\circ {\Phi^{-1}}_{\mid j_B(M(B))}\circ j_B)^{**}(p_B)
\ \geq\ p_A$$
and we have $\Psi^{**}(p_A) = p_B$.
Consequently, $\Psi(\her_{M(A)}(p_A)) = \her_{M(B)}(p_B)$ as required.

\smnoind (b) If $a\in M(A)_{sa}$ and $U$ is an open subset of $\sigma(a) = \sigma(\Phi(j_A(a)))$, then
$\chi_U(\Phi(j_A(a))) = \Phi(\chi_U(j_A(a)))$ is an element of
$\op(B)$ (by \cite[Theorem 2.2]{APT73} and the hypothesis). Thus, by
\cite[Theorem 2.2]{APT73} again, we have $\Phi(j_A(a))\in
j_B(M(B))$. A similar argument shows that
$\Phi^{-1}(j_B(M(B)))\subseteq j_A(M(A))$. Now, we can apply part
(a) to obtain the required conclusion.
\end{prf}

\medskip

\begin{rem}
Note that if $A$ and $B$ are separable and $\Psi: M(A)\to M(B)$ is a $^*$-isomorphism, then $\Psi(A) = B$, by a result of Brown in \cite{Brown88}.
However, the same result is not true if one of them is not separable (e.g.\ take $A = M(B)$ and $\Psi = \id$, where $B$ is non-unital).
Proposition \ref{prop:pre-op-proj}(a) shows that one has $\Psi(A) =B$ if (and only if) $\Psi$ extends to a $^*$-isomorphism from $A^{**}$ to $B^{**}$.
\end{rem}

\medskip

We now consider a suitable equivalence relation on $\op(A)$. A naive
choice is to use the original ``Murray-von Neumann equivalence''
$\sim_{\Mv}$. However, this choice is not good because
\cite{LinHX90} tells us that two open projections that are
Murray-von Neumann equivalent might define non-isomorphic hereditary
$C^*$-subalgebras. On the other hand, one might define
$p\sim_{\her}q$ ($p,q\in \op(A)$) whenever $\her(p) \cong \her(q)$
as $C^*$-algebras. The problem of this choice is that two distinct
open projections of $C([0,1])$ can be equivalent (if they correspond
to homeomorphic open subsets of $[0,1]$), which means that the
resulting classification,  even if possible, will be very different
from the Murray-von Neumann classification.

\medskip

After some thoughts, we end up with an equivalence relation
$\sim_\s$ on $\op(A)$: $p\sim_\s q$ if there
is a partial isometry $v\in A^{**}$ satisfying
\begin{align*}%\label{eq:sp-def-short}
v^*\her_A(p)v\ =\ \her_A(q) \quad \text{and} \quad v\her_A(q)v^*\ =\
\her_A(p).
\end{align*}
Note that this relation is precisely the ``hereditarily stable
version'' of the Murray-von Neumann  equivalence (see Proposition
\ref{prop:cp-equiv}(a)(5) below and the discussion following it).
% and it is stronger than both $\sim_{\Mv}$ and $\sim_{\her}$.

\medskip

In \cite[Definition 1.1]{PZ00}, Peligrad and
Zsid\'{o} introduced another equivalence relation on $\proj(A^{**})$:
$p \sim_{\PZ}q$ if there is a partial isometry $v\in A^{**}$ such that
\begin{align}\label{eqt:PZ}
p=vv^*,\quad q=v^*v,\quad v^*\her_A(p)\subseteq A\quad
\text{and} \quad v\her_A(q)\subseteq A.
\end{align}
It is not difficult to see that $\sim_\PZ$ is stronger than
$\sim_\s$, and a natural description of $\sim_{\PZ}$ on the set of range projections
of positive elements of $A$ is given in \cite[Proposition 4.3]{ORT}.
Moreover, we also gave in \cite[Proposition 3.1]{NW-compare} an equivalent description of $\sim_{\PZ}$ that is similar to $\sim_\s$ but use right ideals instead of hereditary $C^*$-subalgebras.
However, it is now known that $\sim_{\PZ}$ and $\sim_\s$ are actually different even for very simple kind of $C^*$-algebras (see \cite[Theorem 5.3]{NW-compare}).
We decide to use $\sim_\s$ as it
seems to be more natural in the way of using open projections (see Proposition
\ref{prop:cp-equiv}(a) below).

\medskip

Let us start with an extension of $\sim_\s$ to the whole of $\proj(A^{**})$.

\medskip

\begin{defn}\label{def:equ-rel-A}
We say that $p,q \in \proj(A^{**})$ are \emph{spatially equivalent with respect to $A$},
denoted by $p\sim_\s q$, if there exists a partial isometry $v \in A^{**}$ satisfying
\begin{equation}\label{eqt:defn-equiv1}
p = vv^*, \quad q = v^*v, \quad
v^*\her_A(p)v = \her_A(q)
\quad \text{and} \quad
v\her_A(q)v^* = \her_A(p).
\end{equation}
In this case, we also say that the hereditary $C^*$-subalgebras $\her_A(p)$ and $\her_A(q)$ are \emph{spatially
isomorphic}.
\end{defn}

\medskip

It might happen that $\her(p) =0$ but $p\neq 0$ and this is why we
need to consider the first two conditions in \eqref{eqt:defn-equiv1}.
%(in order to avoid a non-zero projection being equivalent to a zero one).
We will see in Proposition
\ref{prop:cp-equiv}(a) that the first two conditions are redundant  if $p$ and $q$ are both open projections.

\medskip

Obviously, $\sim_\s$ is stronger than $\sim_{\Mv}$ (for elements in $\proj(A^{**})$).
Moreover, if $p\sim_\s q$, then $x\mapsto v^*xv$ is a
$^*$-isomorphism from $\her(p)$ to $\her(q)$,  which means that
$\sim_\s$ is stronger than $\sim_{\her}$ in the context of open projections.

\medskip

A good point of the spatial equivalence is
that open projections are stable under $\sim_\s$, as can be seen in part (b) of the following lemma.

\medskip

\begin{lem}\label{lem:sim-prop-equiv}
(a) $\sim_\s$ is an equivalence relation in $\proj(A^{**})$.

\smnoind
(b) Let $p,q\in \proj(A^{**})$ and $u\in A^{**}$ be a partial isometry.
If $p$ is open, $u^*pu = q$, $\her_A(p)\subseteq u\her_A(q)u^*$ and $\her_A(q)\subseteq u^*\her_A(p)u$, then $q$ is open and $p\sim_\s q$.
Consequently, if $p\sim_\s q$ and $p$ is open, then $q$ is open.

\smnoind
(c) If $B\subseteq A$ is a hereditary $C^*$-subalgebra and
$p,q\in \proj(B^{**})$, then $p$ and $q$  are spatially equivalent with respect to $B$ if and only if they are spatially equivalent
with respect to $A$.
\end{lem}
\begin{prf}
(a) It suffices to verify the transitivity. Suppose that
$p,q$ and $v$ are as in Definition \ref{def:equ-rel-A}. If $w\in
A^{**}$ and $r\in \proj(A^{**})$ satisfy that
$$p = w^*w,\ \ r = ww^*,\ \ w\her_A(p)w^* = \her_A(r)\ \ \text{and} \ \ \ w^*\her_A(r)w = \her_A(p),$$
then the
partial isometry $wv$ gives the equivalence $r\sim_\s q$.

\smnoind
(b) As $p$ is open and $\her_A(p)$ is contained in the weak-*-closed subspace
$uA^{**}u^*$, one has $p\leq uu^*$.
Let $v := pu$.
Then $vv^* = p$ and $v^*v =u^*pu = q$.
Moreover, it is clear that $\her_A(p)\subseteq v\her_A(q)v^*$ and $\her_A(q)\subseteq v^*\her_A(p)v$.
Now, it is easy to see that the relations in \eqref{eqt:defn-equiv1} are satisfied.
Furthermore, if $\{a_i\}_{i\in \mathfrak{I}}$ is an approximate unit in $\her_A(p)$, then $\{v^*a_iv\}$ is an increasing net in
$\her_A(q)$ that weak-*-converges
to $v^*pv = q$, and so $q$ is open.
The second statement follows directly from the first one.

\smnoind
(c) Suppose that $p$ and $q$ are spatially equivalent with respect to $A$ and $v\in A^{**}$ satisfies the relations in \eqref{eqt:defn-equiv1}.
As $vv^*, v^*v\in B^{**}$, Remark \ref{rem:her-subalg}(c) tells us that $v\in B^{**}$.
Now the equivalence follows from Remark \ref{rem:her-subalg}(a).
\end{prf}

\medskip

\begin{prop}
\label{prop:cp-equiv}
(a) If $p,q\in \op(A)$, the following statements are equivalent.
\begin{enumerate}[\hspace{1em} (1)]
\item $p\sim_\s q$.
\item $\her(q) = u^*\her(p)u$ and $\her(p) = u\her(q)u^*$ for a partial isometry $u\in A^{**}$.
\item $\her(q) \subseteq u^*\her(p)u$ and $\her(p) \subseteq u\her(q)u^*$ for a partial isometry $u\in A^{**}$.
\item $q\leq v^*v$ and  $v\her(q)v^* = \her(p)$ for a partial isometry $v\in A^{**}$.
\item There is a partial isometry $w\in A^{**}$ such that $p = ww^*$ and
$$\{w^*rw: r\in \op(A); r\leq p\}\ =\ \{s\in \op(A): s\leq q\}.$$ \end{enumerate}

\smnoind (b) If $M$ is a von Neumann algebra and $p,q\in \proj(M)$,
then $p\sim_\s q$ if and only if $p\sim_{\Mv} q$ as elements in $\proj(M)$.
\end{prop}
\begin{prf}
(a) The implications $(1)\Rightarrow (2) \Rightarrow (3)$ and $(1)\Rightarrow (4)$ are clear.

\smnoind $(3)\Rightarrow (1)$.
Since $q$ is open, one has $q\leq u^*u$.
Thus, $(uq)^*uq = q$ and Statement (3) also holds when $u$
is replaced by $uq$.
As $p$ is also open, a similar argument shows
that $p\leq uqu^*$ and Statement (3) holds if we replace $u$ by $v:=puq$ and that $p =
vv^*$.
Furthermore, since $vqv^* = vv^* = p$,  Lemma
\ref{lem:sim-prop-equiv}(b) tells us that $p\sim_\s q$.

\smnoind $(4)\Rightarrow (2)$.
This follows from $v^*\her(p)v = v^*v\her(q)v^*v = \her(q)$.

\smnoind $(1)\Rightarrow(5)$.
Notice that $\op(\her(p)) = \{r\in \op(A):
r\leq p\}$ (see Remark \ref{rem:her-subalg}(b)).
Suppose that $v\in A^{**}$ satisfies \eqref{eqt:defn-equiv1} and $r\in \op(\her(p))$.
If $\{a_i\}_{i\in \mathfrak{I}}$ is an increasing net in $\her(p)$ that $\sigma(A^{**}, A^*)$-converge to $r$, then $\{v^*a_iv\}_{i\in \mathfrak{I}}$ is an increasing net in $\her(q)$ that $\sigma(A^{**}, A^*)$-converge to $v^*rv$ and hence $v^*rv\in \op(\her(q))$.
The argument for the other inclusion is similar.

\smnoind
$(5)\Rightarrow(1)$.
By Statement (5), we have $q =
w^*pw$, and the map $\Phi:x\mapsto w^*xw$  is a $^*$-isomorphism from
$\her(p)^{**}$ to $\her(q)^{**}$.
By Proposition
\ref{prop:pre-op-proj}(b), we see that $\Phi(\her(p)) = \her(q)$ and Statement (4) holds.

\smnoind (b) If $p\sim_\s q$, then $p\sim_\Mv q$ as elements in $\proj(M^{**})$, which implies that $p\sim_\Mv q$ as elements in $\proj(M)$ (by considering the canonical $^*$-homomorphism $\Lambda_M: M^{**}\to M$).
Conversely, if $v\in M$ satisfying $p = vv^*$ and $q=v^*v$, then clearly $v^*\her(p)v = \her(q)$.
\end{prf}

\medskip

One can reformulate Statement (5) of Proposition \ref{prop:cp-equiv}(a) in the following way.

\begin{quotation}
There is a partial isometry $w\in A^{**}$ that induces Murray-von Neumann equivalences between open subprojections of $p$ (including $p$) and open subprojections of $q$ (including $q$).
\end{quotation}

\noindent
Therefore, one may regard $\sim_\s$ as the ``hereditarily stable version'' of the Murray-von Neumann equivalence.
Moreover, if $v\in A^{**}$ satisfies the relations in \eqref{eqt:defn-equiv1}, then by Lemma \ref{lem:sim-prop-equiv}(b), $r\sim_\s v^*rv$ for all $r\in \op(\her(p))$, which means that spatial equivalence is automatically ``hereditarily stable''.
%In fact, if $s := v^*rv$ and $x\in \her(r) \subseteq \her(p)$, then $y:= v^*xv\in \her(q)$ and $sys = v^*rxrv = y$, which means that $y\in \her(s)$. Thus, $\her(r)\subseteq v\her(s)v^*$ and by symmetry $\her(s)\subseteq v^*\her(r)v$.

\medskip

\begin{rem}\label{rem:not-order}
(a) Let $p,q\in \proj(A^{**})$.
We call the unique $p_{\rm int}\in \op(A)$ with $\her(p)=\her(p_{\rm int})$ the \emph{interior} of $p$.
By the bijective correspondence between hereditary $C^*$-subalgebras and open projections, $p_{\rm int}$ is the largest open projection dominated by $p$.
As a direct consequence of Proposition \ref{prop:cp-equiv}(a), we
know that $p_{\rm int} \sim_\s q_{\rm int}$ if and only if
\begin{quotation}
$\her(q) \subseteq u^*\her(p)u$  and $\her(p) \subseteq u\her(q)u^*$ for a partial isometry $u\in A^{**}$.
\end{quotation}

\smnoind
(b) Suppose that $p,q\in \op(A)$.
One might attempt to  define $p\lesssim_{\,\s} q$ if there is $q_1\in \op(A)$
with $p\sim_\s q_1\leq q$.
However, unlike the Murray-von Neumann
equivalence situation, $p\lesssim_{\,\s} q$ and $q\lesssim_{\,\s} p$
does not imply that $p \sim_\s q$.
This can be shown by using a result of Lin.
More precisely, it was shown in
\cite[Theorem 9]{LinHX90} that there exist a separable unital simple
$C^*$-algebra $A$ as well as $p\in \proj(A)$ and $u\in A$ such that
$uu^* =1$ and $p_1=u^*u \leq p$, but $\her(p)$ and
$A$ are not $^*$-isomorphic. In particular, $p\nsim_\s
1$.
Now, we clearly have $p\lesssim_{\,\s} 1$.  On the other hand, as $u\in A$, we
have
$$
u^*Au\ =\ \her(p_1)
\quad \text{and} \quad
u\her(p_1)u^*\ =\ A,
$$
which implies that $1\lesssim_{\,\s} p$.

This example also shows that the same problematic situation appears
even if we replace $\sim_\s$ with the stronger equivalence relation
$\sim_\PZ$ as defined in \eqref{eqt:PZ} (because $u\in A$).
Nevertheless, it was shown in \cite[Theorem 1.13]{PZ00} that a
weaker conclusion holds if one adds an extra assumption on either
$p$ or $q$, but we will not recall the details here.
\end{rem}

\medskip

Let us end this
section with the following well-known example. We give an explicit
argument here for future reference. Note that parts (a) and (b) of it mean
that if $a,b\in A_+$ are equivalent in the sense of Blackadar (i.e.,  there exists $x\in A$ with
$a= x^*x$ and $b=xx^*$; see, e.g., \cite[Definition 2.1]{ORT}),
then their support projections are spatially equivalence (which is also a corollary of \cite[Proposition 4.3]{ORT}, since
$\sim_\PZ$ is stronger than $\sim_\s$).

\medskip

\begin{eg}\label{eg:herx-hera}
Suppose that $x\in A$ with $\|x\| = 1$.
Set $a= x^*x$ and $b=xx^*$.
Let $x = ua^{1/2}$ be the polar decomposition.

\smnoind
(a) It is easy to see that $\overline{aAa} = u^*(\overline{xAx^*})u$ and $\overline{xAx^*} = u(\overline{aAa})u^*$,
%(note that $\overline{aAa} = \overline{a^{1/2}Aa^{1/2}}$)
i.e.,  $\overline{xAx^*}$ is spatially isomorphic to $\overline{aAa}$ (by Proposition \ref{prop:cp-equiv}(a)).

\smnoind
(b) Notice that
$u(\overline{aAa})u^* = \overline{xAx^*}
 \supseteq  \overline{xx^*Axx^*}
 \supseteq  \overline{xx^*xAx^*xx^*}
 \supseteq  u\overline{a^{3/2}Aa^{3/2}}u^*
 =  u(\overline{aAa})u^*$,
and we have $\overline{xAx^*} = \overline{bAb}$.
Similarly, $\overline{x^*Ax}=\overline{aAa}$ and $\overline{x^*A^{**}x}=\overline{aA^{**}a}$, which implies that $\her(x) = \her(a)$.
On the other hand, as $\overline{aAa}$ is a hereditary $C^*$-subalgebra of
$\her(a)$ and $\{a^{1/k}\}_{k\in \mathbb{N}}$ is a sequence in
$\overline{aAa}$ which is an approximate unit for $\her(a)$, one has $\overline{aAa} = \her(a)$.
Consequently, $\her(x) = \overline{x^*Ax}$.

\smnoind
(c) Suppose that $B\subseteq A$ is a hereditary $C^*$-subalgebra and $x\in B$.
Since $\overline{aAa} = \overline{a^2Aa^2}$, we see that $\overline{aBa} = \overline{aAa}$.
Therefore, $\her_B(x) = \her_A(x)$ by part (b).
\end{eg}

%\medskip

%Thus, if $a\in A_+$, then $\her(a)$ is the smallest hereditary $C^*$-subalgebra containing $a$ but in general, $x\notin \her(x)$ (even when $x$ is a partial isometry).

\medskip

\section{$C^*$-semi-finiteness and three types of $C^*$-algebras}

\medskip

As in the case of von Neumann algebras (\cite{Murray-vonNeumann36}), in order to define different ``types''
of $C^*$-algebras, we need to define ``abelian'' and ``finite'' open
projections. ``Abelian'' open projections are defined in the same
way as that of von Neumann algebras. However, in order to define
``finite'' open projections, we need to use our ``hereditarily stable version''  of
Murray-von Neumann equivalence in Section 2.
Note that one cannot
go very far with the original Murray-von Neumann equivalence, because there exist $p,q\in \op(A)$ with $p\sim_\Mv q$
but $\her(p)$ and $\her(q)$ are not isomorphic (see \cite{LinHX90}).
Moreover, one cannot use a direct verbatim translation of the
Murray-von Neumann finiteness.

\medskip

\begin{defn}\label{defn:abel-fin}
(a) Let $q\in \op(A)$ and $p\in \proj(qA^{**}q)$.
The \emph{closure of $p$ in $q$}, denoted by $\bar p^q$,  is the smallest closed projection of $\her(q)$ that dominates $p$.

\smnoind
(b) Let $p,q\in \op(A)$ with $p\leq q$.
The projection $p$ is said to be
\begin{enumerate}[i.]
\item \emph{dense in $q$} if $\bar p^q = q$;
\item \emph{abelian} if $\her(p)$ is a commutative $C^*$-algebra;
\item \emph{$C^*$-finite} if for any $r,s\in \op(\her(p))$ with $r\leq s$ and $r\sim_\s s$, one has $\bar r^s =s$.
\end{enumerate}
If $p$ is dense in $q$, we say that $\her(p)$ is \emph{essential}
in $\her(q)$. We denote by $\op_\CC(A)$ and $\op_\CF(A)$ the set of
all abelian open projections and the set of all $C^*$-finite open
projections of $A$, respectively.
\end{defn}

\medskip

The terminology ``$p$ is dense in $q$'' is used in many places (e.g.\
\cite{PZ00}), while the terminology  ``essential'' comes from
\cite{Zhang89}.

\medskip

Some people might wonder why we do not use the
finiteness as defined in \cite{CP79}. The reason is that we want to
give a classification scheme for $C^*$-algebras using open
projections (and the definition of finiteness in \cite{CP79} seems
not related to open projections). 

\medskip

\begin{rem}\label{rem:defn-C-st-fin}
Let $p\in \op(A)$.

\smnoind (a) Suppose that $p$ is abelian.
If $r,s\in \op(\her(p))$
satisfying $r\leq s$ and $r\sim_\s s$, then $r =s$.
Thus, $p$ is $C^*$-finite.

\smnoind
(b) If $\her(p)$ is finite dimensional, then $p$ is $C^*$-finite.

\smnoind (c) One might ask why we do not define $C^*$-finiteness of
$p$ in the following way: for any $r\in \op(\her(p))$ with $r\sim_\s
p$, one has $\bar r^p = p$. The reason is that the stronger
condition in Definition \ref{defn:abel-fin}(b) can ensure every open
subprojection of a $C^*$-finite projection being $C^*$-finite. Such
a phenomena is automatic for von Neumann algebras.

\smnoind
(d) A hereditary $C^*$-subalgebra $B\subseteq A$ is essential in $A$ if and only if for any non-zero hereditary $C^*$-subalgebra $C\subseteq A$, one has $B\cdot C\neq \{0\}$.
Thus, a closed ideal $I\subseteq A$ is essential in the sense of Definition \ref{defn:abel-fin} if and only it is essential in the usual sense (i.e.,  any non-zero closed ideal of $A$ intersects $I$ non-trivially).
\end{rem}

\medskip

\begin{defn}
A $C^*$-algebra $A$ is said to be:

\begin{enumerate}[i.]
\item \emph{$C^*$-finite} if $1\in\op_\CF(A)$;
\item \emph{$C^*$-semi-finite} if every element in $\op(A)\setminus\{0\}$ dominates an element in $\op_\CF(A)\setminus \{0\}$;
\item \emph{of Type $\ta$} if every element in $\op(A)\cap \Z(A^{**})\setminus\{0\}$ dominates an element in $\op_\CC(A)\setminus \{0\}$;
\item \emph{of Type $\tb$} if $\op_\CC(A) = \{0\}$ but each element in $\op(A)\cap \Z(A^{**})\setminus\{0\}$ dominates an element in $\op_\CF(A)\setminus \{0\}$;
\item \emph{of Type $\tc$} if $\op_\CF(A) = \{0\}$.
\end{enumerate}
\end{defn}

\medskip

Let us give an equivalent form of the above abstract
definition through the relation between open projections (respectively, central open projections) and hereditary $C^*$-subalgebras (respectively, ideals).
%Note that these relations play very important roles in the discussion in this paper.
A $C^*$-algebra $A$
is
\begin{itemize}
\item $C^*$-finite if and only if for each hereditary $C^*$-subalgebra $B\subseteq A$, every hereditary $C^*$-subalgebra of $B$ that is spatially isomorphic to $B$ is essential in $B$;
\item $C^*$-semi-finite if and only if every non-zero hereditary $C^*$-subalgebra of $A$ contains a non-zero $C^*$-finite hereditary $C^*$-subalgebra;
\item of type $\ta$ if and only if every non-zero closed ideal of $A$ contains a non-zero abelian hereditary $C^*$-subalgebra;
\item of type $\tb$ if and only if $A$ does not contain any non-zero abelian hereditary $C^*$-subalgebra and every non-zero closed ideal of $A$ contains a non-zero $C^*$-finite hereditary $C^*$-subalgebra;
\item of type $\tc$ if and only if $A$ does not contain any non-zero $C^*$-finite hereditary $C^*$-subalgebra.
\end{itemize}

\medskip

\begin{rem}\label{rem:simple-type}
Suppose that $A$ is simple.

\smnoind
(a) $A$ is either of type $\ta$, type $\tb$ or type $\tc$.

\smnoind (b) We will see in Corollary \ref{cor:type-postlim} that
$A$ is of type $\ta$ if and only if $A$ is of type I (see, e.g.,
\cite[6.1.1]{Ped79} for its definition). Moreover, if $A$ is of type
II (in the sense of \cite{CP79}), then $A$ is of type $\tb$ (by
Proposition \ref{prop:tr-st>type-B} below), while if $A$ is purely
infinite (in the sense of \cite{Cuntz81}), then $A$ is of type $\tc$
(by Proposition \ref{prop:pure-inf-C}(a) below and \cite[Theorem 1.2(ii)]{Zhang92}). However, we do not
know if the converse of the last two statements hold.
\end{rem}

\medskip

A positive element $a\in A_+$ is said to be \emph{$C^*$-finite} if
$\her(a)$ (i.e.,  $\overline{aAa}$) is $C^*$-finite.
%Parts (a) and (b) of the following results follow from the argument of \cite[Proposition 6.1.7]{Ped79}, but since the settings are slightly different, we give a brief account here.

\medskip

\begin{prop}\label{prop:finite-ideal}
(a) The sum, $\CC(A)$, of all abelian hereditary $C^*$-subalgebras of $A$ is a (not necessarily closed) ideal of $A$.
If $\CC(A)_+ := \CC(A)\cap A_+$, then $\CC(A)$ coincides with the vector space ${\rm span} \ \! \CC(A)_+$ generated by $\CC(A)_+$.

\smnoind
(b) The sum, $\CF(A)$, of all $C^*$-finite hereditary $C^*$-subalgebras of $A$ is a (not necessarily closed) ideal of $A$.
If $\CF(A)_+ := \CF(A)\cap A_+$, then $\CF(A) = {\rm span} \ \! \CF(A)_+$.

\smnoind
(c) If $B\subseteq A$ is a hereditary $C^*$-subalgebra, then $\CC(B)_+=\CC(A)\cap B_+$ and $\CF(B)_+ = \CF(A)\cap B_+$.
\end{prop}
\begin{prf}
Since parts (a) and (b) follow from the arguments of \cite[Proposition 6.1.7]{Ped79}, we will only give the proof for part (c).
Moreover, we will only establish the second equality as the argument for the first one is similar.
As $K_A$ is a hereditary cone, the argument of part (b) tells us that $\CF(A)_+ = K_A$.
It is clear that $\CF(B)\subseteq \CF(A)\cap B$.
Conversely, if $w \in K_A\cap B$ and $w_1,...,w_n\in F_A$ such that $w = \sum_{i=1}^n w_i$, then $w_i\leq w \in B_+$, which implies that $w_i\in F_A\cap B = F_B$ (see Example \ref{eg:herx-hera}(c)).
Consequently, $w\in K_B$ as required.
\end{prf}

\medskip

Clearly, $\CC(A)\subseteq \CF(A)$.
We will see in Theorem \ref{thm:decomp}(d) below that the closed ideal $\overline{\CC(A)}$ is of type $\ta$, while $\overline{\CF(A)}$ is $C^*$-semi-finite.

\medskip

\begin{eg}\label{eg:type-A-finite}
(a) If $A$ is commutative, then $A$ is of type $\ta$ and is $C^*$-finite.
Moreover, $\CC(A) = \CF(A) = A$.

\smnoind
(b) Let $p\in \op(\CL(\ell^2))\subseteq \CL(\ell^2)^{**}$ such that $\her(p) = \CK(\ell^2)$ (the $C^*$-algebra of all compact operators).
Then $p\neq 1$ but $\her(1-p) = (0)$.
In fact, if $T\in \her(1-p)$, we have $pT = 0$ and $ST = SpT = 0$ for any $S\in \mathcal{K}(\ell^2)$, which gives $T=0$.
Moreover, $p$ is dense in $1$ because $\CK(\ell^2)$ is an essential closed ideal of $\CL(\ell^2)$ (see Remark \ref{rem:defn-C-st-fin}(d)).

\smnoind (c) If $H$ is an infinite dimensional Hilbert space, then
$\CK(H)$ is a $C^*$-algebra of type $\ta$, which is not $C^*$-finite
but is $C^*$-semi-finite. In fact, as $\CK(H)$ is simple and
contains many rank-one projections, it is of type $\ta$. On the
other hand, suppose that $e\in \proj(\CK(H))$ is a rank-one
projection. Then $1-e\in \op(\CK(H))\subseteq \CL(H)$ and there is an isometry $v\in
\CL(H)$ with $vv^* = 1-e$. Thus,
$$v^*\her(1-e)v\ =\ \CK(H) \quad \text{and}\quad  1-e\ \sim_\s\ 1.$$
Moreover, as $e\in \proj(\CK(H))$, we see that $1-e$ is also a
closed projection and hence it is not dense in $1$. Finally, as all
hereditary $C^*$-subalgebras  of $\CK(H)$ are given by projections
in $\CL(H)$, they are of the form $\CK(K)$ for some subspaces
$K\subseteq H$. Hence, $\CK(H)$ is $C^*$-semi-finite (see Remark
\ref{rem:defn-C-st-fin}(b)).

\smnoind
(d) Let $H$ be a Hilbert space. Clearly,
$\proj(\CK(H))\subseteq \op_\CF(\CL(H))$.
Hence, if $\mathfrak{F}(H)$ is
the set of all finite rank operators, then $\mathfrak{F}(H)\subseteq
\CF(\CL(H))$. Suppose that $B\subseteq \CL(H)$ is a $C^*$-finite
hereditary $C^*$-subalgebra and $p\in \proj(B)$. As $p$ is
$C^*$-finite and $pBp= p\CL(H)p \cong \CL(K)$ for a subspace
$K\subseteq H$, we see that $K$ is finite dimensional (see part (c)) and so $p\in
\CK(H)$.
Since $B\subseteq \CL(H)$ is a hereditary $C^*$-subalgebra, $B$
is generated by its projections.
Thus, $B$ is a hereditary
$C^*$-subalgebra of $\CK(H)$, and  $B\cong \CK(H')$ for a subspace
$H'\subseteq H$.
The $C^*$-finiteness of $B$ again implies that $\dim H' <
\infty$, and $B\subseteq \mathfrak{F}(H)$. Consequently,
$$\CF(\CL(H))\ =\ \mathfrak{F}(H).$$
On the other hand, since any finite rank projection is a sum of
rank-one projections and any rank-one projection belongs to
$\CC(\CL(H))$, we see that $\mathfrak{F}(H) = \CC(\CL(H)) = \CF(\CL(H))$.
Furthermore, by Proposition
\ref{prop:finite-ideal}(c), we also have $\CF(\CK(H)) = \CC(\CK(H))
= \mathfrak{F}(H)$.
\end{eg}

\medskip

\begin{rem}\label{rem:open-hered-ideal}
Let $e\in \op(A)$ and $z(e)$ be the central
support of $e$ in $A^{**}$.

\smnoind
(a) $z(e) =
\sup_{u\in U_{M(A)}} ueu^*$ (see, e.g., \cite[Lemma 2.6.3]{Ped79}), and
$z(e)$ is an open projection (see Remark \ref{rem:her-subalg}(d)) with $\her(z(e))$
being the smallest closed ideal containing $\her(e)$.

\smnoind
(b) Recall that $B:=\her(e)\subseteq A$ is said to
be \emph{full} if $\her(z(e)) = A$.
In this case, $B$ is strongly Morita
equivalent to $A$ (see, e.g., \cite{Rief}).
Consequently,
$\her(e)$ is always strongly Morita equivalent to $\her(z(e))$.
\end{rem}

\medskip

The following provides an important tool to us in this paper.
An essential ingredient of its proof (in particular, part (b)) is a result of Peligrad and Zsid\'{o} in \cite{PZ00}.

\medskip

\begin{prop}\label{prop:mor-contains}
Let $A$ and $B$ be two strongly Morita equivalent $C^*$-algebras.

\smnoind
(a) $A$ contains a non-zero abelian hereditary $C^*$-subalgebra if and only if $B$ does.

\smnoind
(b) $A$ contains a non-zero $C^*$-finite hereditary $C^*$-subalgebra if and only if $B$ does.
\end{prop}
\begin{prf}
There exist a $C^*$-algebra $D$
and $e\in \proj(M(D))$ such that both $A$ and $B$ are full
hereditary $C^*$-subalgebras of $D$ and we have
$$A \cong eDe \quad \text{and} \quad
B \cong (1-e)D(1-e)$$
(see, e.g., \cite[Theorem II.7.6.9]{Blac}).
Thus, $z(e) = 1 = z(1-e)$.

\smnoind (a) It suffices to show that $A$ contains a non-zero
abelian hereditary $C^*$-subalgebra whenever $D$ does. Let $p\in
\op_\CC(D)\setminus\{0\}$. As $p z(e) = p\neq 0$, we see that
$pueu^* \neq 0$ for some $u\in U_{M(D)}$. By
replacing $p$ with $u^*pu$, we may assume that $pe\neq 0$, and hence
$e\her_D(p)e \neq (0)$. If $x,y\in \her_D(p)$ and $\{b_j\}_{j\in
\mathfrak{I}}$ is an approximate unit of $\her_D(p)$, then $b_i e b_j
\in \her_D(p)$ which implies that
$$xey \ = \ \lim xb_i e b_j y
\ = \ \lim yb_i e b_j x \ = \ yex.$$
Consequently, $e\her_D(p)e$ is an abelian hereditary $C^*$-subalgebra of $A$.

\smnoind (b) It suffices to show that if $D$ contains a non-zero
$C^*$-finite hereditary $C^*$-subalgebra, then so does $A$. Suppose
that $p\in \op_\CF(D)\setminus\{0\}$. By \cite[Theorem 1.9]{PZ00},
there exist $e_0,e_1\in \op(\her_D(e))$ and $p_0,p_1\in \op(\her_D(p))$
satisfying
$$\overline{e_0 + e_1}^e =e, \ \ \overline{p_0 +
p_1}^p =p, \ \ z(e_0)z(p_0) = 0 \ \text{ and }\ e_1 \sim_\PZ p_1.$$
Suppose that $p_1 =0$. Then $e_1 = 0$ and $z(e_0)$ is
dense in $z(e) = 1$ (by \cite[Lemma 1.8]{PZ00}). This implies that
$z(p_0) = 0$, and we have a contradiction that $p_0=0$ is dense in
the non-zero open projection $p$. Therefore, $p_1 \neq 0$ and is
$C^*$-finite.
Since $\her_D(e_1)\cong \her_D(p_1)$ (note that $\sim_\PZ$ is stronger than $\sim_\s$), we see that
$\her_D(e_1)$ is a non-zero $C^*$-finite hereditary $C^*$-subalgebra
of $A = \her_D(e)$.
\end{prf}

\medskip

One may also use the argument of part (b) to obtain part (a), but we keep the alternative argument since it is also interesting.

\medskip

Suppose that $E$ is a full Hilbert $A$-module implementing the
strong Morita equivalence between $A$ and $B$, i.e.,  $B\cong
\CK_A(E)$ (see, e.g., \cite{Lan}).
If $I$ is a closed ideal of $A$, then $EI$ is a full
Hilbert $I$-module and $\CK_I(EI)$ is a closed ideal of $B$.
%This gives an order preserving bijection from the lattice of closed ideals of $A$ to that of $B$.

\medskip

We recall from \cite[Definition 2.1]{PZ00} that $A$
is said to be \emph{discrete}  if any non-zero open projection of
$A$ dominates a non-zero abelian open projection.

\medskip

\begin{thm}\label{thm:morita-types}
(a) Let $A$ and $B$ be two strongly Morita equivalent $C^*$-algebras.
Then $A$ is of type $\ta$ (respectively, type $\tb$ or type $\tc$) if and only if $B$ is of the same type.

\smnoind
(b) A $C^*$-algebra $A$ is of type $\ta$ if and only if it is discrete.
\end{thm}
\begin{prf}
(a) Suppose that $A$ is of type $\tb$. If $\op_\CC(B) \neq \{0\}$,
then $\op_\CC(A)\neq \{0\}$ (because of Proposition
\ref{prop:mor-contains}(a)), which is a contradiction. Let $J$ be a non-zero
closed ideal of $B$. As in the paragraph above, the strong Morita
equivalence of $A$ and $B$ gives a closed ideal $J_0$ of $A$ that is strongly
Morita equivalent to $J$. As $J_0$ contains a non-zero $C^*$-finite
hereditary $C^*$-subalgebra, so is $J$ (by Proposition
\ref{prop:mor-contains}(b)). This shows that $B$ is of type $\tb$.
The argument for the other two types are similar and easier.

\smnoind (b) It suffices to show that if $A$ is of type $\ta$, then
it is discrete. Let $B\subseteq A$ be a non-zero hereditary
$C^*$-subalgebra and $J\subseteq A$ be the closed ideal generated
by $B$ (which is strongly Morita equivalent to $B$; see Remark \ref{rem:open-hered-ideal}(b)). As $J$ contains
a non-zero abelian hereditary $C^*$-subalgebra, so does $B$ (by
Proposition \ref{prop:mor-contains}(a)).
\end{prf}

\medskip

The following result follows from Proposition \ref{prop:mor-contains}(b) and the argument of Theorem \ref{thm:morita-types}.

\medskip

\begin{cor}\label{cor:cp-semi-fin}
(a) $A$ is $C^*$-semi-finite if and only if any non-zero closed ideal of $A$ contains a non-zero $C^*$-finite hereditary $C^*$-subalgebra.

\smnoind
(b) If $A$ is strongly Morita equivalent to a $C^*$-semi-finite $C^*$-algebra, then $A$ is also $C^*$-semi-finite.

\smnoind
(c) $A$ is of type $\tb$ if and only if it is $C^*$-semi-finite and anti-liminary (i.e., it does not contain any non-zero commutative hereditary $C^*$-subalgebra).
\end{cor}

\medskip

\begin{rem}\label{rem:type-hered-ideal}
(a) As in the case of von Neumann algebra, strong Morita equivalence does not preserve $C^*$-finiteness.
In fact, for any $C^*$-algebra $A$, the algebra $A\otimes \CK(\ell^2)$ is not $C^*$-finite (using the same argument as Example \ref{eg:type-A-finite}(c); note that $1\otimes (1-e)$ is both an open and a closed projection of $A\otimes \CK(\ell^2)$).
Consequently, any stable $C^*$-algebra is not $C^*$-finite.

\smnoind
(b) By Remark \ref{rem:open-hered-ideal}(b), Theorem \ref{thm:morita-types}(a) and Corollary \ref{cor:cp-semi-fin}(b), any type $\ta$, type $\tb$, type $\tc$ or $C^*$-semi-finite hereditary $C^*$-subalgebra is contained in a closed ideal of the same type.
\end{rem}

\medskip

Recall that a $C^*$-algebra $A$ has \emph{real rank zero} in the sense of Brown and Pedersen if the set of elements in $A_{sa}$ with finite spectrum is norm dense in $A_{sa}$ (see, e.g., \cite[Corollary 2.6]{BP91}).
The following result follows from Theorem \ref{thm:morita-types}(b), Corollary \ref{cor:cp-semi-fin}(c) as well as the fact that any hereditary $C^*$-subalgebra of a real rank zero $C^*$-algebra is again of real rank zero (see, e.g., \cite[Corollary 2.8]{BP91}).

\medskip

\begin{cor}\label{cor:RR0}
Let $A$ be a $C^*$-algebra with real rank zero.

\smnoind
(a) $A$ is of type $\ta$ if and only if every projection in $\proj(A)\setminus \{0\}$ dominates an abelian projection in $\proj(A)\setminus \{0\}$.

\smnoind
(b) $A$ is of type $\tb$ if and only if every projection in $\proj(A)\setminus \{0\}$ is non-abelian but dominates a $C^*$-finite projection in $\proj(A)\setminus \{0\}$.

\smnoind
(c) $A$ is of type $\tc$ if and only if $A$ does not contain any non-zero $C^*$-finite projection.

\smnoind
(d) $A$ is $C^*$-semi-finite if and only if every projection in $\proj(A)\setminus \{0\}$ dominates a $C^*$-finite projection in $\proj(A)\setminus \{0\}$.
\end{cor}

\medskip

\begin{rem}\label{rem:rr0}
Suppose that $A$ is a $C^*$-finite $C^*$-algebra with real rank zero.
If $r,p\in \proj(A)$ such that $r\leq p$ and there exists $u\in A$ with $uu^*=r$ and $u^*u= p$, then $r\sim_\s p$ and so, $r = \bar r^p = p$.
\end{rem}

\medskip

\begin{cor}\label{cor:CF-RR0}
If $A$ is of real rank zero, then the closures of the ideals $\CC(A)$ and $\CF(A)$ (see Proposition \ref{prop:finite-ideal}) are the closed linear spans of abelian projections and of $C^*$-finite projections in $\proj(A)$, respectively.
\end{cor}
\begin{prf}
If $B\subseteq A$ is a $C^*$-finite hereditary $C^*$-subalgebra, then $B$ is the closed linear span of $\proj(B)\cap \op_\CF(B)$.
Thus, $\CF(A)$ lies inside the closed linear span of $\proj(A)\cap \op_\CF(A)$.
Conversely, it is clear that $\proj(A)\cap \op_\CF(A) \subseteq \CF(A)$.
The argument for the statement concerning $\CC(A)$ is similar.
\end{prf}

\medskip

\begin{cor}\label{cor:type-hered}
Let $A$ be of type $\ta$ (respectively, of type
$\tb$, of type $\tc$ or $C^*$-semi-finite).

\smnoind (a) If $B$ is a hereditary $C^*$-subalgebra of $A$, then
$B$ is of type $\ta$ (respectively, of type $\tb$, of type $\tc$ or $C^*$-semi-finite).

\smnoind (b) If $A$ is a hereditary $C^*$-subalgebra of
$A_0$ that generates an essential ideal $I\subseteq A_0$, then $A_0$ is of type $\ta$ (respectively, of type $\tb$, of type
$\tc$ or $C^*$-semi-finite).
\end{cor}
\begin{prf}
(a) As any hereditary $C^*$-subalgebra of $B$ is a hereditary $C^*$-subalgebra of $A$, this result follows directly from the definitions, Theorem \ref{thm:morita-types}(b) and Corollary \ref{cor:cp-semi-fin}(c).

\smnoind
(b) Note that $A$ is strongly Morita equivalent to $I$ and any hereditary $C^*$-subalgebra of $A_0$ intersects $I$ non-trivially.
Thus, this part follows from the definitions, Theorem \ref{thm:morita-types} and Corollary \ref{cor:cp-semi-fin}.
\end{prf}

\medskip

Consequently, we have the following result.

\medskip

\begin{cor}\label{cor:mult-unit}
Suppose that $A$ is non-unital, and $\ti A$ is the unitalization of $A$.
Then $A$ is of type $\ta$ (respectively, of type $\tb$, of type
$\tc$ or $C^*$-semi-finite) if and only if $\ti A$ is of type
$\ta$ (respectively, of type $\tb$, of type $\tc$ or
$C^*$-semi-finite). The same is true when $\ti A$ is replaced by
$M(A)$.
\end{cor}

\medskip

Our next lemma is probably well-known, but we give a simple argument here for completeness.

\medskip

\begin{lem}\label{lem:proj-her}
Let $e,f\in \op(A)$ and $p,q\in \op(A)\cap \Z(A^{**})$.

\smnoind
(a) $ep\in \op(A)$ and $\her(ep) = \her(e)\cap \her(p)$.

\smnoind
(b) If $e\neq 0$ and $\her(e)\subseteq \her(p) +\her(q)$, then  $\her(e)\cap \her(p) \neq (0)$ or $\her(e)\cap \her(q) \neq (0)$.

\smnoind
(c) If $z(e)z(f) = 0$, then $\her(e) + \her(f) = \her(e+f)$.
\end{lem}
\begin{prf}
Parts (a) and (c) are obvious (see Remark \ref{rem:her-subalg}(d)).
To show part (b), note that as $\her(p) + \her(q)\subseteq \her(p+q-pq)$, we have $e\leq p+q - pq$.
If $ep = 0 = eq$, one obtains a contradiction that $e = e(p+q-pq) = 0$.
Thus, the conclusion follows from part (a).
\end{prf}

\medskip

\begin{lem}\label{lem:sum-C*-fin}
If $\{p_i\}_{i\in \mathfrak{I}}$ is a family in $\op_\CF(A)$ with $z(p_i)z(p_j) = 0$ for $i\neq j$, then $p:=\sum_{i\in \mathfrak{I}} p_i \in \op_\CF(A)$.
\end{lem}
\begin{prf}
It is clear that $p$ is an open projection and $z(p) = \sum_{i\in \mathfrak{I}} z(p_i)$.
Suppose that $r, q\in \op(\her(p))$ with $r\leq q$ and $r\sim_\s q$.
Let $u\in A^{**}$ with $q=u^*u$ and $u\her(q)u^* = \her(r)$.
For any $i\in \mathfrak{I}$, we set $q_i := z(p_i)q, r_i := z(p_i)r\in \op(A)$ and $u_i := z(p_i)u$.
It is easy to see that  $q = \sum_{i\in \mathfrak{I}} q_i$, $r= \sum_{i\in \mathfrak{I}}r_i$, $q_i=u_i^*u_i$ and
$r_i \leq q_i \leq z(p_i)p = p_i$.
By Lemma \ref{lem:proj-her}(c), we see that
$$z(p_i)\her(q)\ =\ z(p_i)\big(\her(q_i) + \her\big({\sum}_{j\in \mathfrak{I}\setminus \{i\}} q_j\big)\big)\ =\ \her(q_i).$$
Similarly, $z(p_i)\her(r) = \her(r_i)$ and we have $u_i\her(q_i)u_i^* = \her(r_i)$.
By Proposition \ref{prop:cp-equiv}(a), we know that $r_i\sim_\s q_i$ and the $C^*$-finiteness of $p_i$ tells us that $r_i$ is dense in $q_i$.
If $e\in \op(\her(q))$ with $re = 0$, then $e_i:= z(p_i)e\in \op(\her(q_i))$ with $r_ie_i = 0$, which means that $e_i = 0$ (because $\overline{r_i}^{q_i} = q_i$).
Consequently, $e = \sum_{i\in \mathfrak{I}} e_i = 0$ and $r$ is dense in $q$ as required.
\end{prf}

\medskip

Part (a) of the following result is the equivalence of statements (i) and (iii) in \cite[Theorem 2.3]{PZ00}, while part (b) follows from the proof of \cite[Theorem 2.3]{PZ00}, Lemma \ref{lem:sum-C*-fin}, Theorem \ref{thm:morita-types}(a) and Corollary \ref{cor:type-hered}(b).

\medskip

\begin{prop}\label{prop:big-finite}
(a) A $C^*$-algebra $A$ is of type $\ta$ if and only if there is an
abelian hereditary $C^*$-subalgebra of $A$ that generates an essential
closed ideal of $A$.

\smnoind (b) A $C^*$-algebra $A$ is $C^*$-semi-finite if and only if
there is a $C^*$-finite hereditary $C^*$-subalgebra of $A$
that generates an essential closed ideal of $A$.
\end{prop}

\medskip

\section{Comparison with existing theories}

\medskip

In this section, we compare our ``Murray-von Neumann type classification'' with existing results in the
literature.
Through these comparisons, we obtain many new
examples of $C^*$-algebras of different types.
Moreover, we will
show that a von Neumann algebra is a
 type $\ta$, type $\tb$, type $\tc$ or
$C^*$-semi-finite $C^*$-algebra if and only if it is, respectively, a type I, type II,
type III or semi-finiteness von Neumann algebra.

\medskip

\subsection{Comparison with type $I$ algebras}\

\bigskip

Recall that a $C^*$-algebra $A$ is said to be of \emph{type I} if for any irreducible representation $(\pi, H)$ of $A$, one has $\CK(H) \subseteq \pi(A)$.
We have already seen in Theorem \ref{thm:morita-types}(b) that type $\ta$ is the same as discreteness.
Thus, the following result is a direct consequence of \cite[Theorem 2.3]{PZ00}.
Note that one can also obtain it using Theorem \ref{thm:morita-types}(a) and \cite[Theorems 1.8 and 2.2]{Beer82}.

\medskip

\begin{cor}\label{cor:postlim>ta}
Any type I $C^*$-algebra is of type $\ta$.
\end{cor}

\medskip

The converse of the above is not true even for real rank zero $C^*$-algebras, as can be seen in the following example.

\medskip

\begin{eg}\label{eg:type-A-not-I}
Example \ref{eg:type-A-finite}(c) and Corollary \ref{cor:type-hered}(b) tell us that $\CL(\ell^2)$ is of type $\ta$.
However, $\CL(\ell^2)$ is not a type I $C^*$-algebra (see, e.g., \cite[6.1.2]{Ped79}).
\end{eg}

\medskip

\begin{prop}\label{prop:I-vs-A}
(a) $A$ is of type I if and only if every primitive quotient of $A$ is of type $\ta$.

\smnoind
(b) If $A$ is of type $\ta$ and contains no essential primitive ideal, then $A$ is of type I.
\end{prop}
\begin{prf}
(a) Because of Corollary \ref{cor:postlim>ta} and the fact that
quotients of type I $C^*$-algebras are also of type I, we only need
to show the ``if'' part.
Let $\pi: A\to \CL(H)$ be an irreducible
representation and $B$ be a non-zero abelian hereditary
$C^*$-subalgebra of $A/\ker\pi$. If $\ti\pi: A/\ker\pi \to \CL(H)$
is the induced representation, the restriction $\ti\pi_B: B \to
\CL(\ti\pi(B)H)$ is non-zero and irreducible. Thus, $\dim \ti\pi(B)H
= 1$ and $\ti\pi(b)$ is a rank-one operator (and hence is
compact) for any $b\in B\setminus \{0\}$. This shows that
$\ti\pi(A/\ker\pi)\cap \mathcal{K}(H)\neq (0)$, and $\pi(A)
\supseteq \mathcal{K}(H)$.

\smnoind (b) Suppose that $\pi: A \to \CL(H)$ is an irreducible
representation and $J$ is a non-zero closed ideal of $A$ with $J
\cap \ker \pi = (0)$. If $B\subseteq J$ is a non-zero abelian
hereditary $C^*$-subalgebra,  the restriction $\pi_B: B \to
\CL(\pi(B)H)$ is non-zero and irreducible. The same argument as in
part (a) tells us that $\pi(A) \supseteq \mathcal{K}(H)$.
\end{prf}

\medskip

\begin{rem}
(a) Proposition \ref{prop:I-vs-A}(a) actually shows that $A$ is of
type I if and only if any primitive quotient contains a  non-zero
abelian hereditary $C^*$-subalgebra, which is likely to be a known fact.

\smnoind (b) If every quotient of $\CL(\ell^2)$ were of type $\ta$,
then Proposition \ref{prop:I-vs-A}(a) told us that $\CL(\ell^2)$ were
a type I $C^*$-algebra, which contradicted \cite[6.1.2]{Ped79}.
Consequently, not every quotient of a type $\ta$ $C^*$-algebra is of type $\ta$.
\end{rem}

\medskip

If $A$ is simple and of type $\ta$, then by Proposition
\ref{prop:I-vs-A}(b), it is of type I. This, together with Example
\ref{eg:type-A-finite}(c), gives the following.

\medskip

\begin{cor}\label{cor:type-postlim}
If $A$ is a simple $C^*$-algebra of type $\ta$, then $A =
\CK(H)$ for some Hilbert space $H$. If, in addition, $A$ is
$C^*$-finite, then $A = M_n$ for some positive integer $n$.
\end{cor}

\medskip

\subsection{Comparison with type $I\!I$ and (semi-)finite $C^*$-algebras}\

\bigskip

The following is a direct consequence of Remark \ref{rem:simple-type}(a) and Corollary \ref{cor:type-postlim}.

\medskip

\begin{cor}\label{cor:inf-dim-C-fin-sim}
Any infinite dimensional $C^*$-finite simple $C^*$-algebra is of type $\tb$.
\end{cor}

\medskip

In the following, we compare type $\tb$ and type $\tc$ with the notions of type $I\!I$ and type $I\!I\!I$ as introduced by Cuntz and Pedersen in \cite{CP79}.
Let us recall from \cite[p.~140]{CP79} that $x\in A_+$ is said to be
\emph{finite} if for any sequence $\{z_k\}_{k\in \mathbb{N}}$ in $A$ with  $x =
\sum_{k=1}^\infty z_k^*z_k$, the condition $\sum_{k=1}^\infty z_kz_k^*\leq
x$ will imply $x = \sum_{k=1}^\infty z_kz_k^*$.
We also recall that $A$ is said to be \emph{finite}
(respectively, \emph{semi-finite}) if every $x\in A_+\setminus \{0\}$ is finite
(respectively, $x$ dominates a non-zero finite element).
Furthermore, $A$ is said to be of \emph{type II} if it is anti-liminary and finite, while $A$ is said to be of \emph{type III} if it has no non-zero finite elements (see \cite[p.~149]{CP79}).

\medskip

Let
$T_s(A)$ be the set of all tracial states on $A$.
It  follows from \cite[Theorem 3.4]{CP79} that $T_s(A)$ separates
points of $A_+$ if $A$ is finite.

\medskip

\begin{prop}\label{prop:tr-st>type-B}
If $T_s(A)$  separates points of
$A_+$, then $A$ is $C^*$-finite. Consequently, if $A$ is finite,
then $A$ is $C^*$-finite.
\end{prop}
\begin{prf}
Suppose on the contrary that there exist $r,q\in \op(A)$ with $r\leq
q$, $r\sim_\s q$ but $\bar r^q \lneq q$. For any $\tau\in T_s(A)$,
if $\ti \tau$ is the normal tracial state on $A^{**}$ extending
$\tau$, then $\ti \tau(r) = \ti\tau(q)$ (because $r = vv^*$ and $q =
v^*v$ for some $v\in A^{**}$). Moreover, if $\{a_i\}_{i\in
\mathfrak{I}}$ is an approximate unit in $\her(r)$, one has $\ti
\tau(r) = \lim \tau(a_i)$. Since $\bar r^q \lneq q$, there exists
$s\in \op(\her(q))\setminus \{0\}$ with $rs = 0$. If $x\in
\her(s)_+$ with $\|x\| = 1$, one can find $\tau_0\in T_s(A)$ with
$\tau_0(x) > 0$.
Thus, we have $\tau_0(a_i) + \tau_0(x) \leq
\ti\tau_0(q)$ (as $a_i + x \leq q$ because $a_i x =0$), which gives the
contradiction that $\ti \tau_0(r) + \tau_0(x) \leq \ti\tau_0(q)$.
\end{prf}

\medskip

As in \cite{CP79}, we denote by $\mathscr{F}^A$ the set of all finite elements in $A_+$.
If $B\subseteq A$ is a hereditary $C^*$-subalgebra, then
$$
\mathscr{F}^B = \mathscr{F}^A\cap B.
$$
In fact, it is obvious that $\mathscr{F}^A\cap B \subseteq \mathscr{F}^B$.
Conversely, suppose that $x\in \mathscr{F}^B$.
Consider $y\in A_+$ and a sequence $\{z_k\}_{k\in \mathbb{N}}$ in $A$ satisfying $y\leq x$, $y = \sum_{k=1}^\infty z_kz_k^*$ and $x = \sum_{k=1}^\infty z_k^*z_k$.
Since $B_+$ is a hereditary cone of $A_+$, we have $y\in B_+$ and $z_k^*z_k, z_kz_k^*\in B_+$ ($k\in \mathbb{N}$).
By Remark \ref{rem:her-subalg}(c), we know that $z_k\in B$ and so, $y =x$ as required.

\medskip

\begin{cor}\label{cor:cp-sf-C-st-sf}
(a) $A$ is semi-finite if and only if every non-zero hereditary $C^*$-subalgebra of $A$ contains a non-zero finite hereditary $C^*$-subalgebra.

\smnoind
(b) If $A$ is semi-finite (respectively, of type II), then $A$ is $C^*$-semi-finite (respectively, of type $\tb$).
\end{cor}
\begin{prf}
(a) For the necessity, let $B\subseteq A$ be a
non-zero hereditary $C^*$-subalgebra.
If $y\in B_+\setminus \{0\}$,
there is $x\in \mathscr{F}^A\setminus \{0\}$ with $x\leq y$.
By \cite[Lemma 4.1]{CP79} and \cite[Theorem
4.8]{CP79} as well as their arguments, one can find a non-zero finite hereditary $C^*$-subalgebra of $\her(x)$.
More precisely, let $f\in C(\sigma(x))_+$ such that $f$ vanishes in a neighborhood of $0$ and $f(t) \leq t \leq f(t) + \frac{\|x\|}{2}$ ($t\in \sigma(x)$).
There exists $g\in C(\sigma(x))_+$ and $\lambda > 0$ such that $f = fg$ and $g(t) < \lambda t$  ($t\in \sigma(x)$).
Then $g(x) \in \mathscr{F}^A$ and $f(x) = f(x)g(x)$, i.e.,
$$f(x)\ \in\ \mathscr{F}_0\ :=\ \{a\in A_+: a = ay \text{ for some } y\in \mathscr{F}^A\}\ \subseteq\ \mathscr{F}^A.$$
For any $z\in \her(f(x))_+$, we have $zg(x) = z$ and $z\in \mathscr{F}_0\cap \her(f(x))\subseteq \mathscr{F}^A\cap \her(f(x)) = \mathscr{F}^{\her(f(x))}$.
Thus, $\her(f(x))$ is a non-zero finite hereditary $C^*$-subalgebra of $\her(x)$.

For the sufficiency, let $y\in A_+\setminus \{0\}$ and $C$ be a
non-zero finite hereditary $C^*$-subalgebra of $\her(y)$. Observe
that $C_+ = \mathscr{F}^C = \mathscr{F}^A\cap C$. Take any $x\in
C_+$ with $\|x\| =1$. Since $x^{1/2}yx^{1/2} \leq \|y\| x\in
\mathscr{F}^A$, we know, from \cite[Lemma 4.1]{CP79}, that
$$y^{1/2}xy^{1/2}\ =\
y^{1/2}x^{1/2}(y^{1/2}x^{1/2})^*
\ \in\ \mathscr{F}^A.$$
Moreover, as $y^{1/2}xy^{1/2}\leq y$, we see that $A$
is semi-finite.

\smnoind
(b) This follows from part (a), Proposition \ref{prop:tr-st>type-B} and Corollary \ref{cor:cp-semi-fin}(c).
\end{prf}

\medskip

\begin{eg}
(a) If $A$ is an infinite dimensional simple $C^*$-algebra with a faithful tracial state, then $A$ is of type $\tb$ (by Corollary \ref{cor:inf-dim-C-fin-sim} and Proposition \ref{prop:tr-st>type-B}).
In particular, if $\Gamma$ is an infinite discrete group such that $C^*_r(\Gamma)$ is simple (see, e.g., \cite{BCH94} for some examples of such groups), then $C^*_r(\Gamma)$ is of type $\tb$.

\smnoind
(b) Every simple $AF$ algebra which is not of the form $\CK(H)$ is of type $\tb$ (because of \cite[Proposition 4.11]{CP79} as well as Corollaries \ref{cor:type-postlim} and \ref{cor:cp-sf-C-st-sf}(b)).
\end{eg}

\medskip

\subsection{Comparison with type $I\!I\!I$ and purely infinite $C^*$-algebras}\

\bigskip

If a $C^*$-algebra $A$ contains a non-zero (positive) finite element $x$, the argument of the necessity of Corollary \ref{cor:cp-sf-C-st-sf}(a) tells us that there is a non-zero finite hereditary $C^*$-subalgebra of $A$, and hence $A$ is not of type $\tc$, because of Proposition \ref{prop:tr-st>type-B}.
This gives the following corollary.

\medskip

\begin{cor}\label{cor:C>III}
If $A$ is of type $\tc$, then it is of type III.
\end{cor}

\medskip

In the following, we will also compare type $\tc$ with the notion of pure infinity as defined by Cuntz (in the case of simple $C^*$-algebras) and by Kirchberg and R\o rdam (in the general case).
Suppose that $a\in M_n(A)$ and $b\in M_m(A)$ ($m,n\in \mathbb{N}$).
As in \cite[Definition 2.1]{KR00}, we say that \emph{$a\precsim b$ relative to $M_{m,n}(A)$} if there is a sequence $\{x_k\}_{k\in \mathbb{N}}$ in $M_{m,n}(A)$ such that $\|x_k^*bx_k - a\| \to 0$.
An element $a\in A$ is said to be \emph{properly infinite} if $a\oplus a\precsim a$ relative to $M_{1,2}(A)$.
Moreover, $A$ is said to be \emph{purely infinite} if every element in $A_+$ is properly infinite (see \cite[Theorem 4.16]{KR00}).
Note that if $A$ is simple, this notion coincides with the one in \cite{Cuntz81}, namely, every hereditary $C^*$-subalgebra of $A$ contains a non-zero infinite projection (see, e.g., the work of Lin and Zhang in \cite{LZ}).

\medskip

\begin{prop}\label{prop:pure-inf-C}
(a) If $A$ has real rank zero and is purely infinite, then it is of type $\tc$.

\smnoind
(b) If $A$ is a separable purely infinite $C^*$-algebra with stable rank one, then $A$ is of type $\tc$.
\end{prop}
\begin{prf}
(a) By \cite[Theorem 4.16]{KR00}, any element $p \in \proj(A)\setminus\{0\}$ is properly infinite and hence is infinite, in the sense that there exist $q\in \proj(A)$ and $v\in A$ such that $q\leq p$, $v^*v = p$ and $q=vv^*$ (see, e.g., \cite[Lemma 3.1]{KR00}).
Thus, $p\sim_\s q$ (as $v\in A$) but $q$ is not dense in $p$ (because $p-q\in \proj(A)\setminus\{0\}$).
Consequently, any non-zero projection in $A$ is not $C^*$-finite, and Corollary \ref{cor:RR0}(c) shows that $A$ is of type $\tc$.

\smnoind
(b) Suppose on contrary that $A$ contains a non-zero $C^*$-finite hereditary $C^*$-subalgebra $B$ and we take any $z\in B_+$ with $\|z\| =1$.
By \cite[Theorem 4.16]{KR00}, one has $z\oplus z \precsim z\oplus 0$ relative to $M_2(A)$, and so, $z\oplus z \precsim z\oplus 0$ relative to $M_2(\her(z))$ (by \cite[Lemma 2.2(iii)]{KR00}).
Thus, \cite[Proposition 4.13]{ORT} implies
$$p_z\oplus p_z
\ =\ p_{z\oplus z}
\ \precsim_{\operatorname{Cu}}\ p_{z\oplus 0}
\ =\ p_z\oplus 0$$
(see \cite[\S 3]{ORT} for the meaning of $\precsim_{\operatorname{Cu}}$).
Moreover, one obviously has $p_{z\oplus 0} \precsim_{\operatorname{Cu}} p_{z\oplus z}$.
Since $A$ has stable rank one, we conclude that
$p_z\oplus p_z \sim_\PZ p_z\oplus 0$
(by \cite[6.2(1)'\&(2)']{ORT}) and hence $p_z\oplus p_z \sim_\s p_z\oplus 0$.
This means that $M_2(\her(z))$ is spatially isomorphic (and hence $^*$-isomorphic) to its hereditary $C^*$-subalgebra $\her(z)\oplus (0)$, which is not essential in $M_2(\her(z))$ (because $(0) \oplus \her(z)$ is a non-zero hereditary $C^*$-subalgebra and we can apply Remark \ref{rem:defn-C-st-fin}(d)).
As $\her(z)$ is $^*$-isomorphic to $\her(z) \oplus (0)$ and hence to $M_2(\her(z))$, we know that $\her(z)$ is also spatially isomorphic to an inessential hereditary $C^*$-subalgebra.
Consequently, $\her(z)$ is not $C^*$-finite, which contradicts the fact that $B$ is $C^*$-finite.
\end{prf}

\medskip

One may regard parts (a) and (b) of the above as two extremes, because any real rank zero $C^*$-algebras has plenty of projections, while a purely infinite $C^*$-algebra with stable rank one is stably projectionless.
Let us make the following conjecture.

\medskip

\begin{conjecture}\label{conj:typeC}
{Every purely infinite $C^*$-algebra is of type
$\tc$.}
\end{conjecture}

\medskip

On the other hand, by Proposition \ref{prop:pure-inf-C} and Corollary \ref{cor:C>III}, we know that any separable purely infinite $C^*$-algebra $A$ having real rank zero or stable rank one is of type III.
This implication actually holds without these extra assumptions, as can be seen in the following proposition, which gives another evidence for Conjecture \ref{conj:typeC}.
Note that this proposition also implies \cite[Proposition 4.4]{KR00}.
To show this result, let us recall the following notation from \cite[p.~3476]{ORT}.
For any $\epsilon > 0$, let $f_\epsilon: \mathbb{R}_+ \to \mathbb{R}_+$ be the function
$$
f_\epsilon(t)\ =\ \begin{cases}
t/\epsilon \ & \text{if } t\in [0,\epsilon)\\
1 & \text{if } t\in [\epsilon, \infty).
\end{cases}
$$
If $\mu\in T_s(A)$ and $a\in A_+$, we define
$$
d_\mu(a)\ :=\ {\sup}_{\epsilon > 0}\ \!\mu(f_\epsilon(a))
$$
(note that the definition in \cite{ORT} is for tracial weights but we only need tracial states here).

\medskip

\begin{prop}\label{prop:pure-inf>III}
Any purely infinite $C^*$-algebra $A$ is of type III.
\end{prop}
\begin{prf}
Suppose on the contrary that $\mathscr{F}^A\neq \{0\}$.
By the argument of the necessity of Corollary \ref{cor:cp-sf-C-st-sf}(a), there is $z\in A_+$ with $\|z\| = 1$ and $\her(z)$ being a finite $C^*$-algebra.
By the argument of Proposition \ref{prop:pure-inf-C}(b), one has $z\oplus z \precsim z\oplus 0$ relative to $M_2(\her(z))$.
By \cite[Remark 2.5]{ORT}, we see that $d_\mu(z\oplus z) \leq d_\mu(z\oplus 0)$ for each $\mu\in T_s(M_2(\her(z)))$.
Now, if $\tau\in T_s(\her(z))$, then $\tau\otimes {\rm Tr_2}\in T_s(M_2(\her(z)))$ (where $\rm Tr_2$ is the canonical tracial state on $M_2$), and the above tells us that
$$ \sup_{\epsilon >0} \tau(f_\epsilon(z))
 =
\sup_{\epsilon >0} (\tau\otimes {\rm Tr_2})(f_\epsilon(z)\oplus f_\epsilon(z))
 \leq  \sup_{\epsilon >0} (\tau\otimes {\rm Tr_2})(f_\epsilon(z)\oplus 0)
 =  \sup_{\epsilon >0} \frac{\tau(f_\epsilon(z))}{2},$$
which gives $d_\tau(z) = 0$ and hence $\tau (z) = 0$.
This contradicts \cite[Theorem 3.4]{CP79}.
\end{prf}

\medskip

If one can show that $\her(a)$ is not $C^*$-finite, for every properly infinite positive element $a$ in any $C^*$-algebra, then the above conjecture is verified.
Let us recall from  \cite[Proposition 3.3(iv)]{KR00} that $a\in A_+$ is properly infinite if and only if there are sequences $\{x_n\}_{n\in \mathbb{N}}$ and $\{y_n\}_{n\in \mathbb{N}}$ in $\her(a)$ such that $x_n^*x_n\to a$, $y_ny_n^*\to a$ and $x_n^*y_n\to 0$.
The following remark tells us that if $a\in A_+$ satisfies a stronger condition than the above, then $\her(a)$ is indeed non-$C^*$-finite.

\medskip

\begin{rem}
Let $a\in A_+$ such that there exist $x,y\in \her(a)$ with $x^*x = a = y^*y$ as well as $x^*y = 0$.
By Example \ref{eg:herx-hera}(a)\&(b), we see that $\her(a)$ is spatially isomorphic to its hereditary $C^*$-subalgebra $\her(x^*)$.
As $\her(x^*)\her(y^*)= (0)$, we see that $\her(x^*)$ is not essential in $\her(a)$.
Thus, $\her(a)$ is not $C^*$-finite.
%If one can show that the same is true for every properly infinite element $a\in A_+$, then by \cite[Theorem 4.16]{KR00}, every purely infinite $C^*$-algebra is of type $\tc$.
\end{rem}

\medskip

\begin{eg}
For any $AF$-algebra $B$, the $C^*$-algebra $\mathcal{O}_2\otimes B$ is purely infinite (by \cite[Proposition 4.5]{KR00}) and is of real rank zero (by \cite[Theorem 3.2]{BP91}), which means that $\mathcal{O}_2\otimes B$ is of type $\tc$ (by Proposition \ref{prop:pure-inf-C}(a)).
Note that one may replace $\mathcal{O}_2$ with any unital, simple, separable, purely infinite, nuclear $C^*$-algebra (which has real rank zero because of \cite[Theorem 1.2(ii)]{Zhang92}).
\end{eg}

\medskip

\subsection{The case of von Neumann algebras}\

\bigskip

In this subsection, we consider the case of von Neumann algebras.
Let us start with the following lemma.
Note that the necessity of part (a) of this result follows directly from Proposition \ref{prop:tr-st>type-B}, but we give an alternative proof here as this argument is also interesting (see Remark \ref{rem:vN-type} below).

\medskip

\begin{lem}\label{lem:abel-fin}
(a) Let $M$ be a von Neumann algebra. Then $p\in \proj(M)$ is finite as
a projection in $M$ if and only if it is $C^*$-finite.

\smnoind
(b) The ideal $\CF(M)$ in Proposition \ref{prop:finite-ideal} is a dense subalgebra of the ideal $J(M)$ generated by finite projections (as defined in \cite{HKNZ}).
\end{lem}
\begin{prf}
(a) Assume that $p$ is finite. Let $\Lambda_M: M^{**}\to M$ be the canonical $^*$-epimorphism.
If $q\in \op(pMp)$, then $\her_M(q) \subseteq \her_M(\Lambda_M(q))$ and $\Lambda_M(q)\leq p$, which imply that $\Lambda_M(q) = \bar q^p$ (notice that $\bar q^p\in pMp$ because of \cite[Theorem II.1]{Akemann70}).

Suppose that $r,q\in \op(pMp)$ such that $r\leq q$ and $r\sim_\s q$.
Consider $w\in M^{**}$ satisfying %$q = ww^*$ and $r = w^*w$.
$$q = ww^*, \quad r = w^*w, \quad w^*\her(q)w = \her(r) \quad
\text{and}\quad  w\her(r)w^* = \her(q).
$$
Define $v:= \Lambda_M(w)$.
Then $\Lambda_M(q) = vv^*$ and $\Lambda_M(r) = v^*v$.
Since $\Lambda_M(r)\leq \Lambda_M(q)\leq p$, the finiteness of $p$ tells us that $\bar r^p = \Lambda_M(r) = \Lambda_M(q) =\bar q^p$.
If $\bar r^q \lneq q$, there is $e\in \op(\her(q))\setminus\{0\}$ with $re = 0$.
Since $e\in \op(\her(p))$, we obtain a contradiction that $\bar r^p \neq \bar q^p$ (as $r\leq p-e$ but $q\nleq p-e$).
This shows that $p$ is $C^*$-finite.

Conversely, if $p$ is $C^*$-finite, then Remark \ref{rem:rr0} implies that $p$ is finite.

\smnoind
(b) This follows from part (a) and Corollary \ref{cor:CF-RR0}.
\end{prf}

\medskip

\begin{rem}\label{rem:vN-type}
(a) Let $p\in M$ be a finite projection.
If $r\in \proj(pMp)$ with $r\sim_\s p$, then Lemma \ref{lem:abel-fin}(a) and Remark \ref{rem:rr0} tell us that $r=p$.
The same is true if we relax the assumption to $r\in \op(pMp)$.
In fact, we first notice that the $C^*$-finiteness of $p$ gives $\bar r^p = p$.
Moreover, suppose that $w\in M^{**}$ and $v\in M$ are as in the proof of Lemma \ref{lem:abel-fin} for the case when $q =p$.
Then $vv^* = p = \bar r^p = v^*v$.
This means that $v$ is a unitary in $pMp$.
As $v\her (r) v^* = \Lambda_M(w\her(r)w^*) =  pMp$, we have $\her(r) = pMp$ and hence $r = p$.

\smnoind
(b) If $A$ is a $C^*$-algebra and $p\in \op(A)$ satisfying $\bar r^p = \bar q^p$ for any $r,q\in \op(\her(p))$ with $r\leq q$ and $r\sim_\s q$, then by the argument of Lemma \ref{lem:abel-fin}, we see that $p$ is $C^*$-finite.
%Note, however, that such condition may not pass onto open subprojection as density is in general not transitive (see \cite[Definition 1.12]{PZ00}).
\end{rem}

\medskip

The following is a direct consequence of Lemma \ref{lem:abel-fin}  and Corollary \ref{cor:RR0}.

\medskip

\begin{thm}\label{thm:cp-types}
Let $M$ be a von Neumann algebra.

\smnoind
(a) $M$ is of type $\ta$ if and only if $M$ is a type I von Neumann algebra.

\smnoind
(b) $M$ is of type $\tb$ if and only if $M$ is a type II von Neumann algebra.

\smnoind
(c) $M$ is of type $\tc$ if and only if $M$ is a type III von Neumann algebra.

\smnoind
(d) $M$ is $C^*$-semi-finite if and only if $M$ is a semi-finite von Neumann algebra.
\end{thm}

\medskip

\section{Factorisations}

\medskip

In this section, we give two factorization type results for general $C^*$-algebras.
Let us first state the following easy lemma.
Notice that if $A$ contains a non-zero abelian hereditary $C^*$-subalgebra $B$, the closed ideal generated by $B$ is of type $\ta$ (by Corollary \ref{cor:type-hered}(b) and Remark \ref{rem:open-hered-ideal}(b)), and the same is true for $C^*$-finite hereditary $C^*$-subalgebra.

\medskip

\begin{lem}\label{lem:ideal}
If $A$ is not of type $\tc$, then $A$ contains  a
non-zero closed ideal of  either type $\ta$ or type $\tb$.
\end{lem}

\medskip

The following is our first factorization type result, which mimics
the corresponding situation for von Neumann algebras.

\medskip

\begin{thm}\label{thm:decomp}
Let $A$ be a $C^*$-algebra.

\smnoind
(a) There is a largest type $\ta$ (respectively, type $\tb$,
type $\tc$ and $C^*$-semi-finite) hereditary $C^*$-subalgebra $J_\ta$
(respectively, $J_\tb$, $J_\tc$ and $J_\SF$) of $A$, which is also an ideal of $A$.

\smnoind
(b) $J_\ta$, $J_\tb$ and $J_\tc$ are mutually disjoint such that
$J_\ta + J_\tb + J_\tc$ is an essential closed ideal of $A$.
If $e_\ta, e_\tb, e_\tc\in \op(A)\cap Z(A^{**})$ with  $J_\ta = \her(e_\ta)$, $J_\tb = \her(e_\tb)$ and $J_\tc = \her(e_\tc)$, then
$$
1 = \overline{e_\ta + e_\tb}^1 + e_\tc.
$$

\smnoind
(c) $J_\ta + J_\tb$ is an essential closed ideal of $J_\SF$.
If $e_\SF \in \op(A)$ with $J_\SF = \her(e_\SF)$, then
$$
e_\SF = \overline{e_\ta}^{e_\SF} + e_\tb.
$$

\smnoind (d) The closure of $\CC(A)$ and $\CF(A)$ (in Proposition \ref{prop:finite-ideal}) are essential
closed ideals of $J_\ta$ and $J_\SF$, respectively.
\end{thm}
\begin{prf}
(a) We first consider the situation of type $\tb$ hereditary $C^*$-subalgebra.
Let $\CJ_\tb$ be the set of all type $\tb$ closed ideals of $A$.
If $\CJ_\tb =\{ (0)\}$, then $J_\tb := (0)$ is the largest type
$\tb$ hereditary $C^*$-subalgebra of $A$ (see Remark \ref{rem:type-hered-ideal}(b)).
Suppose that there exist distinct elements $J_1$ and $J_2$ in $\CJ_\tb$. If $J_1+J_2$ contains a non-zero
abelian hereditary $C^*$-algebra $B$, then by Lemma
\ref{lem:proj-her}(b), one of the two abelian hereditary
$C^*$-subalgebras $B\cap J_1$ and $B\cap J_2$ is non-zero, which
contradicts $J_1, J_2\in \CJ_\tb$.
On the other hand, consider a non-zero closed ideal
$I$ of $J_1+J_2$. Again, by Lemma \ref{lem:proj-her}(b), we may
assume that the closed ideal $I\cap J_1$ is non-zero. Thus,
$I\cap J_1$ contains a non-zero $C^*$-finite hereditary
$C^*$-subalgebra $B$. This shows that $J_1 + J_2\in \CJ_\tb$ and
$\CJ_\tb$ is a directed set.

For any ideal $J$ of $A$, we consider $e_J\in \op(A)\cap \Z(A^{**})$ with $J = \her(e_J)$.
Set
$$J_\tb\ :=\ \overline{{\sum}_{J\in
\CJ_\tb} J}.$$
Then $e_{J_\tb} = w^*$-$\lim_{J\in \CJ_\tb}
e_J$.
If there is $p\in \op_\CC(A)\setminus \{0\}$ such that $\her(p)\subseteq
J_\tb$, then
$$
p\ =\ p e_{J_\tb}\ =\ p e_{J_\tb} p\ =\ w^*\text{-}{\lim}_{J\in \CJ_\tb} p e_J p,
$$
and one can find $J\in \CJ_\tb$ with the abelian algebra
$\her(p)\cap J$ being non-zero (because of Lemma
\ref{lem:proj-her}(a)), which is absurd.
On the other hand, suppose that $I$ is a non-zero closed ideal of $J_\tb$. The argument above tells us that $I\cap J \neq (0)$ for some $J\in \CJ_\tb$, and hence it contains a
non-zero $C^*$-finite hereditary $C^*$-subalgebra. Consequently,
$J_\tb\in \CJ_\tb$.
Finally, if $B\subseteq A$
is a hereditary $C^*$-subalgebra of type $\tb$, then, by Remark \ref{rem:type-hered-ideal}(b), one has $B\subseteq J_\tb$.

The arguments for the statements concerning $J_\ta$, $J_\tc$ and $J_\SF$
are similar and easier.

\smnoind (b) The first statement follows directly from Lemma
\ref{lem:ideal} (any non-type $\tc$ ideal interests either $J_\ta$ or $J_\tb$).
For the second statement, one obviously has
$e_\ta + e_\tb \leq 1 - e_\tc$.
Suppose that $p\in \op(A)$ with
$e_\ta + e_\tb \leq 1 -p$.
We have $p(e_\ta + e_\tb) = 0$.
If $p\nleq e_\tc$, then $\her(p)$ will contain a hereditary $C^*$-subalgebra of either type $\ta$ or type $\tb$ (by Lemma \ref{lem:ideal}) and Lemma \ref{lem:proj-her}(a) will give a contradiction that either $pe_\ta \neq 0$ or $pe_\tb \neq 0$.
Thus, $1 - e_\tc$ is the smallest closed
projection dominating $e_\ta + e_\tb$.

\smnoind
(c) This follows from a similar (but easier) argument as part (b).

\smnoind (d) Clearly, $\CF(A)
\subseteq J_\SF$ and $\CC(A)\subseteq J_\ta$ (see Remark \ref{rem:type-hered-ideal}(b)).
Their closure are both essential because of Proposition \ref{prop:big-finite}.
\end{prf}

\medskip

By Proposition \ref{prop:big-finite}, there is an abelian (respectively, a $C^*$-finite) hereditary $C^*$-subalgebra that generates an essential ideal of $J_\ta$ (respectively, of $J_\tb$).
Moreover, by \cite[Theorem
2.3(vi)]{PZ00}, the largest type $I$ closed ideal $A_{postlim}$ of
$A$ is an essential ideal of $J_\ta$.

\medskip

\begin{rem}\label{rem:comp-id}
For any closed ideal $J$ of $A$, we write $J^\bot$ for the closed ideal $\{a\in A: aJ = (0)\}$.
It is easy to see that if $J_0$ is an essential ideal of $J$, then $J_0^\bot = J^\bot$.

\smnoind
(a) $J_\ta^\bot = A_{postlim}^\bot$ is the largest anti-liminary hereditary $C^*$-subalgebra of $A$ (note that $\overline{aJ_\ta a}$ is a hereditary $C^*$-subalgebra of $J_\ta$ for every $a\in A_+$).
Furthermore, $J_\tb + J_\tc$ is an essential ideal of $J_\ta^\bot$ (by Lemma \ref{lem:ideal}).

\smnoind
(b) $J_\SF^\bot = (J_\ta + J_\tb)^\bot = J_\tc$.

\smnoind (c) $J_\ta^\bot \cap J_\SF = J_\tb$ (compare with Corollary \ref{cor:cp-semi-fin}(c)).
\end{rem}

\medskip

From now on, we denote by $J^A_\ta$, $J^A_\tb$,  $J^A_\tc$ and
$J^A_\SF$, respectively,  the largest type $\ta$, the largest type
$\tb$, the largest type $\tc$ and the largest $C^*$-semi-finite
closed ideals of a $C^*$-algebra $A$.

\medskip

The following is a direct application of Theorem \ref{thm:cp-types}.

\medskip

\begin{cor}\label{cor:decomp}
Let $M$ be a von Neumann algebra.
If $M_I$, $M_{\text{II}}$ and $M_{\text{III}}$ are respectively the type I summand, the type II summand and the type III summand of $M$, then $J^M_\ta = M_I$, $J^M_\tb = M_{\text{II}}$ and $J^M_\tc = M_{\text{III}}$.
\end{cor}

\medskip

Our next theorem is the second factorization type result, which seems to be more interesting for $C^*$-algebra (c.f.\ \cite[Proposition 4.13]{CP79}).

\medskip

\begin{thm}\label{thm:ideal-decomp}
Let $A$ be a $C^*$-algebra.

\smnoind
(a) $A/J^A_\tc$ is $C^*$-semi-finite and $A/(J^A_\ta)^\bot$ is of type $\ta$.

\smnoind
(b) If $A$ is $C^*$-semi-finite, then $A/J^A_\tb$ is of type $\ta$.
\end{thm}
\begin{prf}
(a) Assume, without loss of generality, that $A/J^A_\tc \neq (0)$ and consider $Q: A \to A/J^A_\tc$ to be the canonical map.
Let $I$ be a non-zero closed ideal of $A/J^A_\tc$ and $J:= Q^{-1}(I)$.
Since $J \supsetneq J^A_\tc$, one knows that $J$ contains a non-zero $C^*$-finite hereditary $C^*$-subalgebra $B$.
Since $B\cap J^A_\tc = (0)$, the $^*$-homomorphism $Q$ restricts to an injection on $B$.
Thus, $Q(B)\subseteq I$ is also a non-zero $C^*$-finite hereditary $C^*$-subalgebra, and $A/J^A_\tc$ is $C^*$-semi-finite (by Corollary \ref{cor:cp-semi-fin}(a)).
The proof of the second statement is similar.

\smnoind
(b) This follows from part (a) and Remark \ref{rem:comp-id}(c).
\end{prf}

\medskip

\begin{rem}
Let $\mathcal{S}$ be a statement concerning $C^*$-algebras that is stable under extensions of $C^*$-algebras (i.e.\ if $I$ is a closed ideal of a $C^*$-algebra $A$ such that $\mathcal{S}$ is true for both $I$ and $A/I$, then $\mathcal{S}$ is true for $A$).

\smnoind
(a) If $\mathcal{S}$ is true for all type $\ta$
and all type $\tb$ $C^*$-algebras,
$\mathcal{S}$ is true for all $C^*$-semi-finite
$C^*$-algebras.
If, in addition, $\mathcal{S}$ is true for all type
$\tc$ $C^*$-algebras, it is true for all
$C^*$-algebras.

\smnoind
(b) If $\mathcal{S}$ is true for all discrete $C^*$-algebras and all anti-liminary $C^*$-algebras, then $\mathcal{S}$ is true for all $C^*$-algebras.
\end{rem}

\medskip

The following results follows from Theorem \ref{thm:morita-types}(a).

\medskip

\begin{cor}
If $A$ and $B$ are strongly Morita equivalent, then the closed ideal of $B$ that corresponds to $J^A_\ta$ (respectively, $J^A_\tb$, $J^A_\tc$ and $J^A_\SF$) under the strong Morita equivalence (see the paragraph preceding Theorem \ref{thm:morita-types}) is precisely $J^B_\ta$ (respectively, $J^B_\tb$, $J^B_\tc$ and $J^B_\SF$).
\end{cor}

\medskip

\begin{rem}
It is natural to ask if the closure $\overline{\CC(\cdot)}$ of
$\CC(\cdot)$ (see Proposition \ref{prop:finite-ideal}) is also
stable under strong Morita equivalence. Unfortunately, it is not the
case. Suppose that $A$ is any type I $C^*$-algebra. Then by
\cite[Theorems 1.8 and 2.2]{Beer82}, there is a commutative
$C^*$-algebra $B$ that is strongly Morita equivalent to $A$. Notice
that $\CC(B) = B$ and $\overline{\CC(A)}$ is of type I$_0$ (by
\cite[Proposition 6.1.7]{Ped79}).
Thus, if $\overline{\CC(\cdot)}$
is stable under strong Morita equivalence, then any type I
$C^*$-algebra $A$ will coincide with $\overline{\CC(A)}$ and hence
is liminary (see, e.g., \cite[Corollary 6.1.6]{Ped79}), which is
absurd.
\end{rem}

\medskip

To end this section, we compare $J^A_*$ with $J^{M(A)}_*$.

\medskip

\begin{prop}
(a) If $B\subseteq A$ is a hereditary $C^*$-subalgebra, then $J_\ta^B = J_\ta^A\cap B$, $J_\tb^B = J_\tb^A\cap B$, $J_\tc^B = J_\tc^A\cap B$ and $J_\SF^B = J_\SF^A\cap B$.

\smnoind
(b) $J^{M(A)}_\ta = \{x\in M(A): xA \subseteq J^A_\ta\}$.
Similar statements hold for $J_\tb$, $J_\tc$ and $J_\SF$.

\smnoind
(c) $J^{M(A)}_\tb = \{x\in M(A): xJ^A_\ta = (0) \text{ and } xA \subseteq J^A_\SF\}$

\smnoind
(d) $J^{M(A)}_\tc = \{x\in M(A): xJ^A_\SF = (0)\} = \{x\in M(A): xJ^A_\ta = (0)  \text{ and } xJ^A_\tb = (0)\}$.
\end{prop}
\begin{prf}
(a) Clearly, $J_\ta^B\subseteq B\cap J_\ta^A$.
Conversely, since $B\cap J_\ta^A$ is a type $\ta$ closed ideal of $B$ (by Corollary \ref{cor:type-hered}(a)), we have $B\cap J_\ta^A \subseteq J_\ta^B$.
The other cases follow from similar arguments.

\smnoind
(b) We will only consider the case of $J_\tb$ (since the other cases follow from similar and easier arguments).
Notice that $J^{M(A)}_\tb\cdot A = J^{M(A)}_\tb\cap A = J^A_\tb$ (by part (a)) and
$$J^{M(A)}_\tb\ \subseteq\ J_0\ :=\ \{x\in M(A): xA \subseteq J^A_\tb\}.$$
Suppose that the closed ideal $J_0\subseteq M(A)$ contains a non-zero abelian hereditary $C^*$-subalgebra $B$.
The abelian hereditary $C^*$-subalgebra $B\cap A = B\cdot A\cdot B$ is contained in $J^A_\tb$ and so, $B\cdot A = (0)$, which contradicts the fact that $A$ is essential in $M(A)$ (see Remark \ref{rem:defn-C-st-fin}(d)).
Furthermore, let $I$ be a non-zero closed ideal of $J_0$.
Then $I\cdot A=I\cap A \neq (0)$ and is a closed ideal of $J^A_\tb$.
Thus, $I\cap A$ contains a non-zero $C^*$-finite hereditary $C^*$-subalgebra.
Consequently, $J_0$ is of type $\tb$ and is a subset of $J^{M(A)}_\tb$.

\smnoind
(c) Obviously, $xJ^A_\ta = (0)$ if and only if $xAJ^{A}_\ta = (0)$.
Thus, this part follows from part (b) and Remark \ref{rem:comp-id}(c).

\smnoind
(d) This part follows from a similar argument as part (c) as well as Remark \ref{rem:comp-id}(b).
\end{prf}

\bigskip

\bibliographystyle{plain}

\end{document}